\documentclass[11pt]{article}
\usepackage{style}
\usepackage{float}

\usetikzlibrary{arrows.meta}
\tikzset{
  agn/.style={circle, draw, fill=white, inner sep=2pt, minimum size=3mm, font=\footnotesize},
  ag arrow/.style={-{Stealth[length=1.3mm,width=1.1mm]}},
  ag dashed arrow/.style={dashed,-{Stealth[length=1.3mm,width=1.1mm]}},
}

\newcommand{\agedge}[2]{%
  \ifnum\numexpr#2-#1\relax>1
    \draw[ag arrow] (n#1) to[bend left=30] (n#2);%
  \else
    \draw[ag arrow] (n#1) -- (n#2);%
  \fi}
\newcommand{\agraph}[2][]{%
  \def\agtmp{#1}%
  \begin{tikzpicture}[baseline={([yshift=-.7ex]n2.center)}, line width=0.4pt]
    \node[agn] (n1) at (0,0)    {1};
    \node[agn] (n2) at (1,0) {2};
    \node[agn] (n3) at (2,0)  {3};
    \foreach \e in {#2} {\expandafter\agedge\e}
  \end{tikzpicture}%
}
\newcommand{\agdashededge}[2]{%
  \ifnum\numexpr#2-#1\relax>1
    \draw[ag dashed arrow] (n#1) to[bend left=30] (n#2);%
  \else
    \draw[ag dashed arrow] (n#1) -- (n#2);%
  \fi}
\newcommand{\agraphsetting}[2]{%
  \begin{tikzpicture}[baseline={([yshift=-.7ex]n2.center)}, line width=0.4pt]
    \node[agn] (n1) at (0,0) {1};
    \node[agn] (n2) at (1,0) {2};
    \node[agn] (n3) at (2,0) {3};
    \foreach \e in {#1} {\expandafter\agedge\e}
    \foreach \e in {#2} {\expandafter\agdashededge\e}
  \end{tikzpicture}%
}

\apptocmd{\thebibliography}{\setlength{\itemsep}{0pt}}{}{}

\let\eqref\cref

\newcommand*\samethanks[1][\value{footnote}]{\footnotemark[#1]}

\hypersetup{linkcolor=PineGreen}

\title{Graph splitting algorithms for inconsistent feasibility problems}
\author{
	Rub\'en Campoy\thanks{Department of Mathematics, University of Alicante, \textsc{Spain}, \href{mailto:ruben.campoy@ua.es}{ruben.campoy@ua.es}, \href{mailto:cesar.lopez@ua.es}{cesar.lopez@ua.es}, \href{mailto:david.torregrosa@ua.es}{david.torregrosa@ua.es}.}
	\and C\'esar L\'opez-Pastor\samethanks[1]
	\and David Torregrosa-Bel\'en\samethanks[1]
}
\date{}


\funding{This research was partially supported by grant PID2025-168246NB-I00 funded by  MICIU/AEI/10.13039/501100011033 and by ERDF/EU, as well as by grant CIGE/2026/80 funded by Conselleria d'Educaci\'o, Cultura i Universitats, Generalitat Valenciana. The second author was additionally supported by grant PREP2022-000118 funded by \seqsplit{MICIU/AEI/10.13039/501100011033} and by ``ESF Investing in your future''.}
\interests{The authors declare that they have no conflict of interest.}
\data{}
\referees{}

\begin{document}

\maketitle

\begin{abstract}
In this work, we study graph splitting algorithms for feasibility problems involving closed and convex sets whose intersection may be empty. We characterize their minimal displacement vector, which is known to be the limit of the differences between successive iterates, in terms of a least-squares problem measuring disagreement between the sets according to the graph structure. We then specialize to closed affine subspaces, proving strong convergence of the shadow sequence to a solution of such a least-squares problem whenever one exists. Our analysis extends classical Douglas--Rachford results on best approximation pairs to graph splitting algorithms.
\end{abstract}

\makepaperdata

\section{Introduction}

Given a real Hilbert space $\Hi$ and a family of nonempty closed and convex sets $C_1,\ldots,C_n\subseteq\Hi$, the \emph{feasibility problem} consists of finding a point in their intersection, that is,
\begin{equation}\label{FP}
    \text{find } x^\ast\in\bigcap_{i=1}^n C_i.
\end{equation}
The family of \emph{projection algorithms} is a classical approach to numerically solving such problems. These methods approximate a solution by iteratively computing projections onto each set separately. One of the most renowned among them is the \emph{Douglas--Rachford (DR) algorithm}. Although it traces back to \cite{DR56}, the method was introduced by Lions and Mercier \cite{LM79} to find a zero of the sum of two maximally monotone operators. When specialized to normal cones, it turns into a projection method admitting a particularly simple formulation. Specifically, given any starting point $v^0\in\Hi$, the DR algorithm iterates via the recurrence
\begin{equation}\label{DRit}
\left\{ \begin{aligned}
x_1^{k+1} &=  P_{C_1}(v^k),\\
x_2^{k+1} &=  P_{C_2}(2x_1^{k+1}-v^k),\\
v^{k+1} &=   v^k - x_1^{k+1} + x_2^{k+1},
\end{aligned}\right.
\end{equation}
where $P_{C}$ denotes the metric projector onto a set $C$. Under the assumption $C_1\cap C_2\neq\emptyset$, the \emph{governing sequence} $(v^k)_{k\geq 0}$ converges weakly to a fixed point of the associated DR operator, whereas the \emph{shadow sequences} $(x_1^k)_{k\geq 1}$ and $(x_2^k)_{k\geq 1}$ converge weakly to the same point in $C_1\cap C_2$, which solves \eqref{FP}.

However, in this paper we are interested in \emph{inconsistent feasibility problems}, where the intersection may be empty. In the case of two sets, a natural replacement when no common feasible point exists is to solve
\begin{equation}\label{eq:best-aprox-pair}
    \min_{(x_1,x_2)\in C_1\times C_2}\frac12\|x_1-x_2\|^2.
\end{equation}
Solutions to the above problem, that is, pairs of points achieving the minimum distance between the sets, are called \emph{best approximation
pairs}. The behavior of DR in this setting was first studied by Bauschke, Combettes, and Luke~\cite{BCL2004}. They showed that the shadow sequences remain bounded and their weak cluster points provide  best approximation pairs, whenever these exist. Subsequently, Bauschke and Moursi~\cite{BM16affine} studied the case of two possibly nonintersecting affine subspaces and established convergence of the entire shadow sequences. This was later generalized to the affine--convex setting by Bauschke, Dao, and Moursi \cite{BDM16affineconvex}. Finally, Bauschke and Moursi~\cite{BM17} culminated the analysis by extending the results to the general convex setting. A central role in these analyses is played by the \emph{minimal displacement vector} of the associated DR operator (see~\cref{def:mdv}). This vector describes the asymptotic difference between successive governing iterates (see~\cref{fact:asymptotic-displacement}) and is closely related to best approximation pairs (see~\cref{rem:DR-bestapprox}).

The first approach for applying DR to feasibility problems involving more than two sets was Pierra's product-space reformulation \cite{pierra}. This reduces the original problem to a two-set feasibility problem in a higher-dimensional Hilbert space, at the cost of increasing the number of governing variables. In the inconsistent case, this approach also provides least-squares approximations; see \cite[\S~3]{BDM16affineconvex}. Other extensions include the \emph{cyclic Douglas--Rachford method}  by Borwein and Tam~\cite{cyclic-dr}, whose behavior for possibly inconsistent problems was studied in~\cite{cyclic-dr-inconsistent}.

More recently, alternative splitting schemes requiring fewer governing variables have been proposed, including the algorithms of Ryu~\cite{ryu20} and Malitsky--Tam~\cite{malitsky2023resolvent}, as well as parallel methods based on a reduced-dimensional product space reformulation~\cite{campoy,condat2023proximal}. A unifying framework for these methods was introduced by Bredies, Chenchene, and Naldi~\cite{graph-drs}, in which the relationships between the variables are represented by a pair of directed graphs.  Under suitable conditions, each pair of graphs defines a convergent scheme, giving rise to the family of \emph{graph splitting algorithms}. This framework not only recovers the aforementioned algorithms but also produces new ones, such as the one generated by the complete graph~\cite{graph-fb,graph-drs}. The benefits of this unified approach are also illustrated in~\cite{graph-linear}, where a more detailed analysis of this family for linear subspaces extends the results of~\cite{bauschke-splitting} without requiring a separate analysis for each algorithm.

Nevertheless, the convergence theory of graph splitting algorithms has been developed under the assumption that the underlying problem admits a solution. In this work, we provide an initial study of this family of methods for the inconsistent case. In particular, we specialize to the setting of feasibility problems described by closed convex sets with possibly empty intersection. We summarize our main contributions as follows:
\begin{itemize}
    \item We characterize the minimal displacement vector and relate it to a least-squares problem that generalizes \eqref{eq:best-aprox-pair} to the context of graph splitting algorithms.
    \item For closed affine subspaces, we establish strong convergence of the full shadow sequence to a solution of the associated least-squares problem, whenever one exists. This addresses the question concerning the inconsistent
affine case raised in \cite[\S~8]{bauschke-splitting}.
\end{itemize}

This paper is organized as follows. \cref{sec:pre} is dedicated to introducing notation and compiling fundamental results. In \cref{sec:convex}, the graph splitting methods are studied for (possibly inconsistent) convex feasibility problems. Particularly, we characterize the displacement vector in terms of a least-squares problem associated with the graph splitting framework. In \cref{sec:affine}, we further investigate the affine case, for which strong convergence results are obtained. We illustrate our results for some instances of graph splitting methods in \cref{sec:illustrative}. Finally, conclusions and open questions are stated in \cref{sec:conc}.

\section{Preliminaries}\label{sec:pre}

Throughout this work, $\Hi$ is a real Hilbert space with inner product $\langle\cdot,\cdot\rangle$ and induced norm $\|\cdot\|$. We use $\to$ and $\rightharpoonup$ to denote strong and weak convergence, respectively.  The Cartesian product Hilbert space $\Hi^n$ is endowed with the inner product
\(
    \langle\bm x,\bm y\rangle
    :=\sum_{i=1}^n\langle x_i,y_i\rangle,
\)
for $\bm x=(x_1,\ldots,x_n)$, $\bm y=(y_1,\ldots,y_n)\in\Hi^n$. The diagonal of $\Hi^n$ is the subspace
\[
    \Delta_n:=\{(x,\ldots,x)\in\Hi^n:x\in\Hi\}.
\]
Elements of $\Hi^n$ and $\mathbb R^n$ will be denoted in bold. The vector of ones in $\R^n$ is denoted by $\bm 1_n$.

\subsection{Operator theory}

For a mapping $T:\Hi\to\Hi$, its \emph{range} and \emph{set of fixed points} are the sets
\[
    \ran T:=\{Tx:x\in\Hi\}
    \quad\text{and}\quad
    \Fix T:=\{x\in\Hi:Tx=x\}.
\]
We denote the closure of $\ran T$ by $\overline{\ran}\,T$. The identity operator is denoted by $\Id$, with its underlying space understood from the context. A mapping $T:\Hi\to\Hi$ is \emph{nonexpansive} if
\[
    \|Tx-Ty\|\leq\|x-y\|,
    \qquad\forall x,y\in\Hi,
\]
and \emph{firmly nonexpansive} if
\[
    \|Tx-Ty\|^2
    +\|(\Id-T)x-(\Id-T)y\|^2
    \leq\|x-y\|^2,
    \qquad\forall x,y\in\Hi.
\]
Every firmly nonexpansive mapping is nonexpansive, and the set of fixed points of a nonexpansive mapping is closed and convex \cite[Corollary~4.24]{bauschke}.

Our analysis concerns mappings which may not have fixed points. For every nonexpansive mapping $T:\Hi\to\Hi$, the set
$\overline{\ran}(\Id-T)$ is nonempty, closed, and convex \cite[Example~20.29 and Corollary~21.14]{bauschke}.
Consequently, the following vector is well-defined.

\begin{definition}[Minimal displacement vector]\label{def:mdv}
Let $T:\Hi\to\Hi$ be nonexpansive.
Its \emph{minimal displacement vector} is
\begin{equation}\label{eq:mdv}
    v_T
    :=\operatorname*{argmin}_{v\in\overline{\ran}(\Id-T)}\|v\|.
\end{equation}
In particular, \(\|v_T\|=\inf_{x\in\Hi}\|x-Tx\|\), and $\Fix T\neq\emptyset$ implies $v_T=0$.
\end{definition}

The following result describes the role of the minimal displacement vector in fixed-point iterations, even when the mapping has no fixed points.

\begin{fact}[{\cite[Corollary~2.3]{baillon1978asymptotic}}]\label{fact:asymptotic-displacement}
Let $T:\Hi\to\Hi$ be firmly nonexpansive and let $\theta\in{]0,2[}$. Let $v^0\in\Hi$ and generate
\[
v^{k+1}=(1-\theta) v^k + \theta T(v^k), \quad \forall k\geq 0.
\]
Then $v^k-Tv^k\to v_T$.
\end{fact}

Let $\mathcal K$ be another real Hilbert space and let $L:\Hi\to\mathcal K$ be a continuous linear operator.
Its kernel is $\ker L:=\{x\in\Hi:Lx=0\}$, and its adjoint is the unique continuous linear operator $L^*:\mathcal K\to\Hi$ satisfying
\[
    \langle Lx,y\rangle
    =\langle x,L^*y\rangle,
    \qquad\forall x\in\Hi,\ \forall y\in\mathcal K.
\]
The adjoint verifies that $(\ker L)^\perp=\overline{\ran}\,L^*$ and $(\ran L)^\perp=\ker L^*$ \cite[Fact 2.25]{bauschke}, where $U^\perp$ denotes the orthogonal complement of a linear subspace $U$.

The \emph{Moore--Penrose inverse} of $L$ is denoted by $L^\dagger$. When $L$ is injective with closed range, it can be computed as
\begin{equation}\label{eq:pseudoinvT}
	L^\dagger=(L^*L)^{-1}L^*.
\end{equation}

Given a matrix $K\in\mathbb R^{m\times n}$, we write $\ot{K}$ for the continuous linear operator $K\otimes\Id:\Hi^n\to\Hi^m$, which acts as
\begin{equation}\label{eq: KronMatrix}
    \ot{K}\bm{x}
    =
    \left(\textstyle\sum_{j=1}^n K_{1,j}x_j,\sum_{j=1}^n K_{2,j}x_j,\ldots,\sum_{j=1}^n K_{m,j}x_j\right)\in\Hi^m,
\end{equation}
for $\bs x=(x_1,\ldots,x_n)\in\Hi^n$. In particular, for
$\bm a=(a_1,\ldots,a_n)\in\mathbb R^n$,
\[
    \ot{\bm a}x=(a_1x,\ldots,a_nx),
    \qquad
    \ot{\bm a}^*\bm x=\sum_{i=1}^n a_ix_i.
\]
This canonical embedding respects matrix multiplication, adjoints and inverses; that is, we have that $\ot{(K_1)}\ot{(K_2)}=\ot{(K_1K_2)}$, $(\ot{K})^*=\ot{K^*}$ and $(\ot{K})^{-1}=\ot{K^{-1}}$, whenever the dimensions are compatible.

\subsection{Projections onto convex sets}

Let $C\subseteq\Hi$ be a nonempty closed and convex set. For every $x\in\Hi$, there exists a unique point in $C$ nearest to $x$ \cite[Theorem 3.16]{bauschke}. The \emph{(metric) projector} $P_C:\Hi\to\Hi$ assigns to each point in $\Hi$ its unique nearest point in $C$, that is,
\[
    P_C(x):=\operatorname*{argmin}_{c\in C}\|x-c\|.
\]
The distance from $x$ to $C$ is \(\operatorname{dist}(x,C):=\inf_{c\in C}\|x-c\|=\|x-P_C(x)\|\). The projector onto a closed and convex set is firmly nonexpansive \cite[Proposition~4.16]{bauschke}. Moreover, projections onto translated sets satisfy \cite[Proposition~3.19]{bauschke}
\begin{equation}\label{eq:traslationPC}
    P_{a+C}(x)=a+P_C(x-a),
    \qquad\forall x,a\in\Hi.
\end{equation}

The \emph{normal cone} of $C$ is defined at $x\in C$ as
\[
    N_C(x):=
        \{u\in\Hi:
          \langle u,c-x\rangle\leq0,\, \forall c\in C\},
\]
whereas $N_C(x)=\emptyset$ when $x\not\in C$. The normal cone is related to the metric projection by \cite[Proposition~6.47]{bauschke}
\begin{equation}\label{fact:proj_nc}
p=P_C(x) \quad\Longleftrightarrow\quad x-p\in N_C(p),\quad \forall x,p\in\Hi.
\end{equation}

Let $C_1,\ldots,C_n\subseteq\Hi$ be nonempty closed convex sets and set $\mathcal C:=C_1\times\cdots\times C_n$. %
For all $\bm x = (x_1,\ldots,x_n)\in\Hi^n$, it holds that $N_\mathcal{C}(\bm x)=N_{C_1}(x_1) \times\cdots \times N_{C_n}(x_n)$  \cite[Proposition~16.9 and Example~16.13]{bauschke}, while $P_{\mathcal{C}}(\bm x) = \left( P_{C_1}(x_1),\ldots, P_{C_n}(x_n)\right)$ \cite[Proposition~29.3]{bauschke}.

When $U\subseteq\Hi$ is a closed linear subspace, the projector $P_U$ is linear, $P_U+P_{U^\perp}=\Id$ and $\ker P_U=U^\perp$ \cite[Corollary~3.24]{bauschke}. Moreover, for a closed affine subspace $C=a+U$,  its normal cone is \(N_C(x)=U^\perp\) for all $x\in C$ \cite[Example~6.43]{bauschke}.

\subsection{Graph theory}

A \emph{directed graph} $G=(\mathcal N, \mathcal E)$ is given by a finite set $\mathcal N$ of \emph{nodes} and a collection of \emph{edges} $\mathcal E$, which are (ordered) pairs of $\mathcal N \times \mathcal N$. The number of nodes $|\mathcal N|$ defines the \emph{order} of the graph. A  \emph{spanning subgraph} of $G$ is a subgraph  whose nodes coincide with the nodes of $G$, namely, $G^\prime = (\mathcal N, \mathcal E^\prime)$ with $\mathcal E^\prime \subseteq \mathcal E$.

Two nodes $i,j\in\mathcal N$ are said to be adjacent if $(i,j)\in\mathcal E$ or $(j,i)\in\mathcal E$.  The \emph{in-degree}, the \emph{out-degree} and the \emph{degree} of a node $i\in\mathcal N$ are respectively defined as 
\[
d_i^{\rm in} = | \{j\in\mathcal N \,:\, (j,i)\in\mathcal E \}|,
\quad d_i^{\rm out} = | \{j\in\mathcal N \,:\, (i,j)\in\mathcal E \}| \quad \text{and} \quad
d_i = d_i^{\rm in} + d_i^{\rm out}.
\]
A graph is \emph{connected} if there exists a \emph{path} between any pair of distinct nodes. By a path we refer to a finite sequence of distinct nodes $\nu=(\nu_1,\ldots,\nu_r)$ with $r\geq 2$ such that $\nu_k$ and $\nu_{k+1}$ are adjacent for all $k=1,\ldots, r-1$. In such a case, we say that $\nu$ is a path between $\nu_1$ and $\nu_r$.

We shall study splitting algorithms that build upon a special class of graphs, ensuring well-defined update orders.

\begin{definition}[Algorithmic graph]
We say that $G=(\mathcal N,\mathcal E)$ is an \emph{algorithmic graph} if (i) $\mathcal N = \{1,\ldots,n\}$ with $n\geq 2$, (ii) $(i,j)\in\mathcal E \Rightarrow i < j$, and (iii) $G$ is connected.
\end{definition}

Let $G=(\mathcal N, \mathcal E)$ be an algorithmic graph of order $n$. The \emph{Laplacian matrix} of  $G$ is the matrix $\operatorname{Lap}(G)\in\R^{n\times n}$ given componentwise by
\begin{equation*}\label{eq:Laplacian}
\operatorname{Lap}(G)_{i,j}:= 
\left\{
\begin{array}{ll}
d_i &\text{if } i=j, \\
-1 &\text{if $i$ and $j$ are adjacent}, \\
0 & \text{otherwise.}
\end{array}
\right.
\end{equation*}
The Laplacian matrix of an algorithmic graph admits an \emph{onto decomposition}, in the sense that
\begin{equation}\label{eq:Zdecom}
\operatorname{Lap}(G) = ZZ^*, \quad \text{with } Z\in\R^{n\times (n-1)},
\end{equation}
which implies $\ker Z=\{0\}$ and $\ker Z^* =  \operatorname{span}\{\mathbf{1}_n\}$, see \cite[Fact~2.17]{graph-linear}. In addition, by \cite[Proposition~2.18]{graph-linear}, the linear and continuous operator $\ot{Z}:\Hi^{n-1}\to\Hi^n$ constructed as in~\cref{eq: KronMatrix} satisfies
\begin{equation}\label{fact:ranZ}
\begin{aligned}
    \ran \ot{Z}
    &=\Delta_n^\perp
    =\left\{\bm w\in\Hi^n: w_1+\cdots+w_n=0\right\},\\
    \ker \ot{Z^*}
    &=\Delta_n
    =\left\{\ot{(\bm 1_n)}u: u\in\Hi\right\}.
\end{aligned}
\end{equation}
In particular, $\ran \ot{Z}$ is closed. Since $\ot{Z}$ is also injective, then $\ot{Z}^\dagger=\ot{(Z^\dagger)}$ can be computed using~\eqref{eq:pseudoinvT}.

\section{Graph splitting methods for convex feasibility}\label{sec:convex}

In this section, we present a family of projection algorithms whose variable-update rules are determined by algorithmic graphs. This family arises  from the graph-based splitting framework developed by Bredies, Chenchene and Naldi in~\cite{graph-drs} specialized to feasibility problems of the form~\cref{FP}. The variable-update rules in the resulting algorithms will be determined by the specification of a \emph{graph setting}, which we introduce next.

\begin{definition}[Graph setting]
\label{def:graph-setting}
Let $G=(\mathcal N,\mathcal E)$ be an algorithmic graph of order $n$, and denote by $d_i$ the degree of node $i\in\mathcal N$. Let $G'=(\mathcal N,\mathcal E')$ be a connected spanning subgraph of~$G$, and fix an onto decomposition
\begin{equation*}
\label{eq:Z-decomposition}
    \operatorname{Lap}(G')=ZZ^*,
    \qquad
    Z\in\mathbb R^{n\times(n-1)},
\end{equation*}
as described in \cref{eq:Zdecom}. We refer to $(G,G',Z)$ as a
\emph{graph setting}.
\end{definition}

\cref{alg:T} displays the graph-based projection method determined by a graph setting. Observe that the $i$-th shadow variable $x_i^{k+1}$ is computed as the projection of a linear combination of shadow and governing variables. Specifically, the shadow variable $x_h^{k+1}$ contributes to the computation of $x_i^{k+1}$ if and only if  $(h,i)\in\mathcal E$, that is, there is an edge of $G$ going from $h$ to $i$. Note that the ordering of the algorithmic graph (condition (ii) in its definition) ensures that the update rule is well-posed. On the other hand, the onto decomposition of $\operatorname{Lap}(G^\prime)$ determines how shadow variables ($x_i^k$) and governing variables ($v_j^k$) affect each other. For a detailed interpretation, we refer the reader to~\cite[Remark~3.8]{graph-fb}.

\begin{algorithm}[ht!]
    \caption{Graph splitting method for convex feasibility problems}\label{alg:T}
    \begin{algorithmic}[1]
    	\State \textbf{input:} a family of nonempty closed convex sets $C_1,\ldots,C_n\subseteq\Hi$, a graph setting
$(G,G',Z)$ as in \cref{def:graph-setting} (with degrees $d_1,\ldots,d_n$), and a relaxation parameter $\theta\in{]0,2[}$.
        \State \textbf{initialize:} $v_1^0,\ldots,v_{n-1}^0\in\Hi$.
        \For {$k=0,1,2,\ldots$}
        \State {$x_i^{k+1}=P_{C_i}\left(\frac{2}{d_i}\sum_{(h,i)\in\mathcal{E}}x_h^{k+1}+\frac{1}{d_i}\sum_{j=1}^{n-1}Z_{i,j}v^k_j\right)$}, \hfill{$i= 1,\ldots,n,$}
        \State {$v_i^{k+1}=v_i^k-\theta\sum_{j=1}^n Z_{j,i}x_j^{k+1}$}, \hfill{$i=1,\ldots,n-1$.}
        \EndFor
    \end{algorithmic}
\end{algorithm}

Let us first express \cref{alg:T} as a relaxed fixed-point iteration. To this end, for a fixed graph setting $(G,G',Z)$ and the family $\bm C=\{C_i\}_{i=1}^n$, consider the mapping \(\mathcal{X}_{\bm C}:\Hi^{n-1}\to\Hi^n\), whose components are recursively defined by
\begin{equation}\label{eq:XC}
    \left(\mathcal X_{\bm C}(\bm v)\right)_i
    :=
    P_{C_i}\left(
        \frac{2}{d_i}
        \sum_{(h,i)\in\mathcal E}\left(\mathcal X_{\bm C}(\bm v)\right)_h
        +\frac{1}{d_i}(\ot{Z} \bm v)_i
    \right),
    \quad \forall i=1,\ldots,n.
\end{equation}
We then define the associated graph operator
\begin{equation}\label{eq:T-graph}
    \T_{\bm C}
    :=
    \Id
    -\ot{Z}^*\mathcal X_{\bm C}.
\end{equation}
With this notation, \cref{alg:T} can be equivalently written as the relaxed fixed-point iteration
\begin{equation}\label{eq:T-iteration}
    \bm v^{k+1}
    =
    (1-\theta)\bm v^k+\theta\T_{\bm C}(\bm v^k), \quad \forall k\geq 0.
\end{equation}

The operator $\T_{\bm C}$ defined in \cref{eq:T-graph} coincides with $\widetilde \T$ in \cite[Eq.~(16)]{graph-linear}, specialized to normal cones (see~\cite[Proposition~3.5]{graph-linear}). Thus, we can exploit the explicit description of its set of fixed points established in that work. To express this characterization, we introduce the \emph{degree balance} of $G$, defined by
\begin{equation}\label{eq:delta}
\bm\delta:=(\delta_1,\ldots,\delta_n)\in\mathbb R^n,
        \quad\text{with } 
        \delta_i:=d_i^{\rm out}-d_i^{\rm in},\quad \forall i=1,\ldots,n,
\end{equation}
as well as the vector
\begin{equation}\label{eq:alpha}
\bm\alpha:=Z^\dagger\bm\delta\in\mathbb R^{n-1}.
\end{equation}
Since $\sum_{i=1}^n\delta_i=0$, we have $\bm\delta\in\{\bm1_n\}^\perp=\ran Z$. In fact, $\bm\alpha$ is the unique vector satisfying $Z\bm\alpha=\bm\delta$ (see~\cite[Proposition~3.9(i)]{graph-linear}).

\begin{lemma}\label{lem:T-firmly}
The operator $\T_{\bm C}$ defined in \cref{eq:T-graph} is firmly nonexpansive and
\begin{equation*}\label{eq:FixT}
    \Fix\T_{\bm C}
    =
    \left\{
        \ot{\bm\alpha}x+\ot{Z}^\dagger\bm u
        \;:\;\textstyle
        x\in\bigcap_{i=1}^n C_i,\quad
        u_i\in N_{C_i}(x), \forall i=1,\ldots,n,\quad
        \sum_{i=1}^n u_i=0
    \right\},
\end{equation*}
where $\bm\alpha\in\R^{n-1}$ is given by \eqref{eq:alpha}.
\end{lemma}
\begin{proof}
Combining \cite[Proposition~3.5]{graph-linear} and \cite[Lemma~2.3]{graph-drs}, we see that $\T_{\bm C}$ is the resolvent of a maximally monotone operator and, thus, firmly nonexpansive by \cite[Corollary 23.9]{bauschke}. The expression for $\Fix\T_{\bm C}$ follows from \cite[Remark~3.10]{graph-linear} applied to $A_i=N_{C_i}$, for $i=1,\ldots,n$.
\end{proof}

As a consequence of \cref{eq:T-iteration}, \cref{lem:T-firmly} and \cref{fact:asymptotic-displacement}, we obtain the following result.

\begin{corollary}[Graph residual convergence]
\label{cor:graph-residual}
Let $(\bm v^k)_{k\geq0}$ and $(\bm x^k)_{k\geq1}$ be the governing and shadow sequences, respectively, generated by \cref{alg:T}. Then
\begin{equation*}\label{eq:residual-convergence}
    \ot{Z}^*\bm x^{k+1}
    \to
    v_{\T_{\bm C}},
\end{equation*}
equivalently, $
    \bm v^k-\bm v^{k+1}
    \to
    \theta v_{\T_{\bm C}}.$
\end{corollary}

\begin{remark}[The consistent case]\label{rem:consistent}
Suppose that the feasibility problem \cref{FP} is consistent; that is, \(\cap_{i=1}^n C_i\neq\varnothing\). Then \cref{lem:T-firmly} gives $\Fix\T_{\bm C}\neq\varnothing$, and hence $v_{\T_{\bm C}}=0$. Thus, \cref{cor:graph-residual} yields
\[
    \ot{Z}^*\bm x^{k+1}\to0.
\]
In particular, the components of the shadow sequence asymptotically reach consensus. Indeed, by \cite[Eq.~(3.25)]{graph-drs},
\[
    \lambda_+
    \operatorname{dist}(\bm x^{k+1},\Delta_n)^2
    \leq
    \|\ot{Z}^*\bm x^{k+1}\|^2,\quad \forall k\geq 0,
\]
where $\lambda_+>0$ is the smallest positive eigenvalue of $\operatorname{Lap}(G')$ (known as the \emph{algebraic connectivity} of $G'$). Hence, \(\operatorname{dist}(\bm x^{k+1},\Delta_n)\to0\).
\end{remark}

In this paper, we focus on the inconsistent case and study the minimal displacement vector~$v_{\T_{\bm C}}$, which determines the limiting residual in \cref{cor:graph-residual}. Our first step is to characterize the range of the displacement mapping $\Id-\T_{\bm C}$ in terms of the product set
\begin{equation*}
    \mathcal C:=C_1\times\cdots\times C_n\subseteq\Hi^n
\end{equation*}
and the structure of the graph setting $(G,G',Z)$. 

\begin{proposition}[Range of the displacement mapping]\label{prop:range-convex}
Let $\T_{\bm C}$ be given by \cref{eq:T-graph}. Then
\begin{equation*}
\ran(\Id-\T_{\bm C})
=
\left\{
\ot{Z}^*\bm x \;:\;
 \bm x=(x_1,\ldots,x_n)\in \mathcal C \;\;\text{with }\;
 \ot{\bm \delta}^*\bm x\in \sum_{i=1}^n N_{C_i}(x_i)
\right\},
\end{equation*}
where $\bm\delta\in\R^{n}$ is given by \eqref{eq:delta}.
\end{proposition}

\begin{proof} By~\eqref{eq:T-graph}, $\ran(\Id-\T_{\bm C}) = \ran \ot{Z}^*\X_{\bm C} = \ot{Z}^*(\ran \X_{\bm C})$. Hence, it  suffices to characterize the range of $\X_{\bm C}$. By \cref{eq:XC} and \cref{fact:proj_nc}, $\bs x=(x_1,\ldots,x_n)\in\ran \X_{\bm C}\subseteq \mathcal C$ if and only if there exist  $\bm v\in\Hi^{n-1}$ and $u_i\in N_{C_i}(x_i)$, for $i=1,\ldots,n$, such that
\begin{equation}\label{eq:ui_NCi}
    (\ot{Z}\bm v)_i
    =
    u_i
    -2\sum_{(h,i)\in\mathcal E}x_h
    +d_i x_i,
    \quad \forall i=1,\ldots,n.
\end{equation}
Define $\bm r=(r_1,\ldots,r_n)\in\Hi^n$ componentwise by \( r_i := u_i-2\sum_{(h,i)\in\mathcal E}x_h+d_i x_i\), for $i=1,\ldots,n$. Then there exists $\bm v$ satisfying \cref{eq:ui_NCi} if and only if $\bm r\in\ran\ot{Z}$, which is equivalent to $\sum_{i=1}^n r_i=0$ in view of \cref{fact:ranZ}. Using $d_i=d_i^{\rm out}+d_i^{\rm in}$ and \eqref{eq:delta}, we have
\begin{align*}
    \sum_{i=1}^n r_i
    &=
     \sum_{i=1}^n u_i
    -2\sum_{i=1}^n\sum_{(h,i)\in\mathcal E}x_h
    +\sum_{i=1}^n
      d_ix_i\\
    &=\sum_{i=1}^n u_i
    -2\sum_{i=1}^n d_i^{\rm out}x_i
    +\sum_{i=1}^n
      (d_i^{\rm out}+d_i^{\rm in})x_i
    =
    \sum_{i=1}^n u_i
    -\sum_{i=1}^n\delta_i x_i.
\end{align*}
Hence, we have proved that
\[
\ran \X_{\bm C}
=
\left\{
 \bm x=(x_1,\ldots,x_n)\in \mathcal C \;:\;
 \ot{\bm \delta}^*\bm x\in \sum_{i=1}^n N_{C_i}(x_i)
\right\},
\]
which concludes the proof.
\end{proof}

\begin{remark}[The Douglas--Rachford case]
\label{rem:DR-range}
The classical Douglas--Rachford method is recovered by taking the unique algorithmic graph of order $2$, with $\mathcal E=\{(1,2)\}$. In this case, $G'=G$, $d_1=d_2=1$, and we may choose $Z=(1,-1)^*$. We obtain $\mathcal X_{\bm C}(v)=\bigl(P_{C_1}(v),P_{C_2}(2P_{C_1}(v)-v)\bigr)$, and hence
\[ \T_{\bm C} = \Id-P_{C_1}+P_{C_2}R_{C_1},\quad \text{with } R_{C_1}:=2P_{C_1}-\Id.\]
Thus, \cref{alg:T} with $\theta=1$ reduces to the classical Douglas--Rachford method \eqref{DRit}. Moreover, the degree balance is $\bm\delta=(1,-1)$. Consequently, \cref{prop:range-convex} yields
\[
    \ran(\Id-\T_{\bm C})
    =
    \{x_1-x_2:x_1\in C_1,\ x_2\in C_2,\
      x_1-x_2\in N_{C_1}(x_1)+N_{C_2}(x_2)\}.
\]
This recovers the normal-cone specialization of the classical Douglas--Rachford displacement-range formula; see \cite[Corollary~2.14]{bauschke2014generalized}.
\end{remark}

\begin{example}[A nonclosed displacement range]
Consider the Douglas--Rachford graph setting described in \cref{rem:DR-range} and the nonempty closed convex sets
\[
    C_1=\mathbb R\times\{0\},
    \qquad
    C_2=\{(s,t)\in\mathbb R^2:t\geq e^s\}.
\]
For each $k\geq 1$ set $z^k:=(-k,-e^{-k})$. Then we have
\[
    P_{C_1}(z^k)=(-k,0),
    \qquad
    P_{C_2}\bigl(R_{C_1}(z^k)\bigr)=(-k,e^{-k}),
\]
and therefore \((\Id-\T_{\bm C})z^k=(0,-e^{-k})\to(0,0)\) as $k\to\infty$. However, $(0,0)\notin\ran(\Id-\T_{\bm C})$, since ${C_1}\cap{C_2}=\varnothing$. Thus, $\ran(\Id-\T_{\bm C})$ is not closed.
\end{example}

When the feasibility problem is inconsistent, we cannot find solutions to \eqref{FP}. A natural alternative is to seek points $x_i\in C_i$, for $i=1,\ldots,n$, that come as close as possible to reaching consensus. As we will show next, within the graph-based framework, this disagreement is measured with respect to the topology of $G'$, by summing the squared distances between points associated with adjacent nodes. Since $ZZ^*=\operatorname{Lap}(G')$, this yields the following least-squares problem:
\begin{equation}\label{eq:graph-ls}\tag{$P_{\mathrm{LS}}$}
    \min_{\bm x\in\mathcal C}
    \frac12\|\ot{Z}^*\bm x\|^2
    =
    \min_{\substack{x_i\in C_i\\i=1,\ldots,n}}
    \frac12\sum_{(i,j)\in\mathcal E'}
    \|x_i-x_j\|^2.
\end{equation}

In what follows, we aim to show that, in fact, the norm of the minimal displacement vector $v_{\T_{\bm C}}$ determines the optimal value for \eqref{eq:graph-ls}.

\begin{lemma}\label{lem:graph-ls}
For every $\bm y\in\mathcal C$, there exists some $v_{\bm y}\in\ran(\Id-\T_{\bm C})$ such that $\norm{v_{\bm y}}\leq\norm{\ot{Z^*}\bm y}$.
\end{lemma}

\begin{proof}
Since $\bm\delta\in\{\bm 1_n\}^\perp = \ran Z = \ran ZZ^* = \ran \Lap(G')$, there exists $\bm\eta^0\in\mathbb R^n$ satisfying $\Lap(G')\bm\eta^0=-\bm\delta$. Set
\[
    \bm\eta:=\bm\eta^0-m\bm 1_n,\quad \text{with } m:=\min_{1\leq i\leq n}\eta_i^0.
\]
Then $\eta_i\geq0$ for all $i=1,\ldots,n$. Moreover, since $\Lap(G')\bm 1_n=0$,
\[
    \Lap(G')\bm\eta
    =
    \Lap(G')\bm\eta^0-m\Lap(G')\bm 1_n
    =
    -\bm\delta.
\]
Also, recall that $\bm\delta\neq0$, so $\bm\eta\notin\ker(\Lap(G'))=\spa\{\bm 1_n\}$, that is,
\begin{equation}\label{cap zero}
    \spa\{\bm 1_n\}\cap\spa\{\bm\eta\}=\{0\}.
\end{equation}
Let $I_n$ be the identity matrix of dimension $n$ and let us denote by $M:= I_n -\bm\eta\bm\eta^*/\|\bm \eta\|^2\in\R^{n\times n}$ the symmetric positive semidefinite matrix defining the orthogonal projection onto the subspace~$\{\bm\eta\}^\perp$ (thus $\ker M=\spa\{\bm \eta\}$, and $M^2=M$), see~\cite[Example~29.18]{bauschke}. The matrix $B\in\R^{n\times n}$, defined as
\[B:=\Lap(G')+ M,\]
is symmetric, positive semidefinite and satisfies that $B\bm\eta=\Lap(G')\bm\eta=-\bm\delta$. More precisely, we can show that $B\succ 0$. Indeed, assume that $\bm z^*B\bm z=0$ for some $\bm z\in\R^n$. Then,
\[ 0 = \bm z^*B\bm z = \bm z^*\Lap(G')\bm z + \bm z^*M\bm z = \|Z^*\bm z\|^2+ \|M\bm z\|^2,\]
so $\bm z\in\ker Z^*\cap\ker M=\spa\{\bm 1_n\}\cap\spa\{\bm\eta\}$. This implies $\bm z=0$ in view of \cref{cap zero}.

Now, fix any $\bm y\in\mathcal C$ and define the convex quadratic function $\Phi_{\bm y}:\Hi^n\to\R$ as
\[\Phi_{\bm y}(\bm x):= \frac12\norm{\ot Z^*\bm x}^2+\frac12\|\ot M(\bm x-\bm y)\|^2=\frac12\langle \ot B \bm x,\bm x\rangle - \langle \ot M \bm y,\bm x\rangle+\frac12\langle \ot M\bm y,\bm y\rangle.\]
Observe that
\( \langle\ot{B}\bm x,\bm x\rangle \geq \lambda_{\min}(B)\|\bm x\|^2\), for all $\bm x\in\Hi^n$,
where $\lambda_{\min}(B)>0$ denotes the smallest eigenvalue of $B$. Then, $\Phi_{\bm y}$ is coercive and \cite[Proposition 11.15]{bauschke} ensures that the function $\Phi_{\bm y}$ has a minimizer $\overline{\bm x}_{\bm y}$ over $\mathcal C$. Moreover, by \cite[Proposition~27.8]{bauschke}, $\overline{\bm x}_{\bm y}\in\mathcal C$ is characterized by the optimality condition
\begin{equation*}\label{eq:optimality_xy}
    \bm u:= -\nabla\Phi_{\bm y}(\overline{\bm x}_{\bm y})=-\ot B\overline{\bm x}_{\bm y}+\ot M\bm y
    \in
    N_{\mathcal C}\left(\overline{\bm x}_{\bm y}\right).
\end{equation*}
Since $\bm\eta^*M=0$, $B=B^*$, and $B\bm\eta=-\bm\delta$, we obtain that
\[
    \ot{\bm\eta}^*\bm u
    =  -\ot{(\bm\eta^*B)}\overline{\bm x}_{\bm y}  +\ot{(\bm\eta^*M)}\bm y
    = -\ot{(B\bm\eta)}^*\overline{\bm x}_{\bm y}
    = \ot{\bm\delta}^*\overline{\bm x}_{\bm y}.
\]
By the nonnegativity of $\bm \eta$, it then follows that
\[
\ot{\bm\delta}^*\overline{\bm x}_{\bm y}
 = \sum_{i=1}^n\eta_i u_i \in \sum_{i=1}^n N_{C_i}((\overline{\bm x}_{\bm y})_i),
\]
which proves that $v_{\bm y}:=\ot{Z}^*\overline{\bm x}_{\bm y}\in\ran(\Id-\T_{\bm C})$ according to \cref{prop:range-convex}. Hence, it holds that 
\[
   \frac12\|v_{\bm y}\|^2= \frac12\|\ot{Z}^*\overline{\bm x}_{\bm y}\|^2\leq\Phi_{\bm y}(\overline{\bm x}_{\bm y})
    \leq\Phi_{\bm y} (\bm y) = \frac12\norm{\ot Z^*\bm y}^2,
\]
proving the desired inequality.
\end{proof}

The following is the main result of this section.

\begin{theorem}[Minimal displacement vector and least-squares problem] \label{thm:graph-ls}
Let $\T_{\bm C}$ be the operator in \cref{eq:T-graph}. Then
\begin{equation}\label{eq:vT-projection}
    v_{\T_{\bm C}}=P_{\overline{\ot{Z}^*(\mathcal C)}}(0),
    \quad\text{ equivalently }\quad
    \inf_{\bm x\in\mathcal C}\frac12\|\ot{Z}^*\bm x\|^2=\frac12\|v_{\T_{\bm C}}\|^2.
\end{equation}
\end{theorem}

\begin{proof}
    To prove \cref{eq:vT-projection}, it suffices to check that $v_{\T_{\bm C}}\in\overline{\ot Z^*(\mathcal C)}$ and that $\norm{v_{\T_{\bm C}}}\leq\norm{\ot Z^* \bm x}$ for all $\bm x\in\mathcal C$. For proving the former, \cref{prop:range-convex} yields the inclusion $\ran(\Id-\mathcal T_{\bm C})\subseteq\ot Z^*(\mathcal C)$. Hence, by \cref{eq:mdv}, we see that $v_{\T_{\bm C}}\in\overline{\ran}(\Id-\mathcal T_{\bm C})\subseteq\overline{\ot Z^*(\mathcal C)}$.

    Now, for any $\bm x\in \mathcal C$, by \cref{lem:graph-ls} there exists some $v_{\bm x}\in\ran(\Id-\T_{\bm C})$ satisfying that $\norm{v_{\bm x}}\leq\norm{\ot Z^*\bm x}$. But by the definition of $v_{\T_{\bm C}}$, we obtain $\norm{v_{\T_{\bm C}}}\leq\norm{v_{\bm x}}\leq\norm{\ot Z^*\bm x}$.
\end{proof}

\begin{remark}[Best approximation pairs]\label{rem:DR-bestapprox}
In the Douglas--Rachford setting of \cref{rem:DR-range}, the least-squares problem \cref{eq:graph-ls} reduces to finding a best approximation pair as in \eqref{eq:best-aprox-pair}. Furthermore, \cref{thm:graph-ls} gives
\[v_{\T_{\bm C}}=P_{\overline{C_1-C_2}}(0),\]
recovering \cite[Theorem~3.4]{BCL2004}.
\end{remark}

We finish this section by establishing some characterizations of the existence of minimizers.

\begin{corollary}[Existence of least-squares solutions]
\label{cor:graph-ls}
The following conditions are equivalent:
\begin{enumerate}
    \item\label{thm:vT-projectioni} \eqref{eq:graph-ls} attains its minimum;
    \item\label{thm:vT-projectionii} $v_{\mathcal T_{\bm C}}\in\ran(\Id-\mathcal T_{\bm C})$;
    \item\label{thm:vT-projectioniii} $v_{\mathcal T_{\bm C}}\in\ot Z^*(\mathcal C)$.
\end{enumerate}
In addition, the next condition implies all the above ones
\begin{enumerate}\setcounter{enumi}{3}
    \item\label{thm:vT-projectioniv} $\ot Z^*(\mathcal C)$ is closed.
\end{enumerate}
Finally, under these assumptions, it holds that
    \begin{equation}\label{argmin ZC}
    \operatorname*{argmin}_{\bm x\in\mathcal C}
    \tfrac12\|\ot{Z}^*\bm x\|^2
    =
    \left\{
        \bm x\in\mathcal C
        \;:\;
        \ot{Z}^*\bm x=v_{\T_{\bm C}}
    \right\}.
    \end{equation}
\end{corollary}

\begin{proof}
\cref{thm:vT-projectioni}$\Rightarrow$\cref{thm:vT-projectionii}. If \eqref{eq:graph-ls} attains a minimum, then there exists some $\overline{\bm x}\in\mathcal C$ such that $\norm{\ot{Z}^*\overline{\bm x}}=\norm{v_{\T_{\bm C}}}$. Applying \cref{lem:graph-ls}, there exists some $v_{\overline{\bm x}}\in\ran(\Id-\T_{\bm C})$ such that $\norm{v_{\overline{\bm x}}}\leq\norm{v_{\T_{\bm C}}}$. Thus, by \cref{eq:mdv}, $v_{\T_{\bm C}}=v_{\overline{\bm x}}\in\ran(\Id-\T_{\bm C})$.

\cref{thm:vT-projectionii}$\Rightarrow$\cref{thm:vT-projectioniii}. It comes from the inclusion $\ran(\Id-\T_{\bm C})\subseteq\ot{Z}^*(\mathcal C)$.

\cref{thm:vT-projectioniii}$\Rightarrow$\cref{thm:vT-projectioni}. Since $v_{\mathcal T_{\bm C}}=\ot Z^*\overline{\bm x}$ for some $\overline{\bm x}\in\mathcal C$, then \cref{eq:vT-projection} yields that $\overline{\bm x}$ is a minimizer.

\cref{thm:vT-projectioniv}$\Rightarrow$\cref{thm:vT-projectioniii}. By \cref{eq:vT-projection}, $v_{\mathcal T_{\bm C}}\in\overline{\ot{Z}^*(\mathcal C)}=\ot{Z}^*(\mathcal C)$.

Finally, the equivalence \cref{thm:vT-projectioni}$\Leftrightarrow$\cref{thm:vT-projectioniii} provides \eqref{argmin ZC}.
\end{proof}

\section{The case of closed affine subspaces}\label{sec:affine}

In this section, we consider the case where the sets are closed affine subspaces of the form
\begin{equation}\label{eq:Ci-affine}
    C_i=\xi_i+U_i,
    \quad \text{for } i=1,\ldots,n,
\end{equation}
with $\xi_i\in\Hi$ and $U_i\subseteq\Hi$ being a closed linear subspace. Following \cite{graph-linear}, we employ the notation
\begin{align*}
    U&:=U_1\cap\cdots\cap U_n\subseteq\Hi,\\
    \mathcal U&:=U_1\times\cdots\times U_n\subseteq\Hi^n,
\end{align*}
and write $\bm \xi:=(\xi_1,\ldots,\xi_n)$, so that $\mathcal C=\bm \xi+\mathcal U$. Define also the mapping $\X_{\bm U}:\Hi^{n-1}\to\Hi^n$ analogously to \cref{eq:XC}, but related to $\bm U=\{U_i\}_{i=1}^n$; that is, for $\bm v\in\Hi^{n-1}$,
\begin{equation}\label{eq:XU}
    \left(\X_{\bm U}(\bm v)\right)_i
    :=
    P_{U_i}\left(
        \frac{2}{d_i}
        \sum_{(h,i)\in\mathcal E}
        \left(\X_{\bm U}(\bm v)\right)_h
        +
        \frac{1}{d_i}(\ot{Z}\bm v)_i
    \right),
    \quad \forall i=1,\ldots,n,
\end{equation}
and consider the (homogeneous) graph operator
\begin{equation}\label{eq:T-U}
    \T_{\bm U}
    :=
    \Id-\ot{Z}^*\X_{\bm U}.
\end{equation}
Since each $P_{U_i}$ is linear and continuous, so are $\X_{\bm U}$ and $\T_{\bm U}$. This linearity allows us to derive the following translation formula.

\begin{lemma}[Affine translation formula]\label{lem:affine-splitting}
Let $\T_{\bm C}$ and $\T_{\bm U}$ be defined in \cref{eq:T-graph,eq:T-U}, respectively. Then
\begin{equation}\label{eq:T-affine-splitting}
    \T_{\bm C}(\bm v)
    =
    \T_{\bm U}(\bm v)
    -
    \ot{Z}^*\X_{\bm C}(0),
    \qquad \forall \bm v\in\Hi^{n-1}.
\end{equation}
\end{lemma}
\begin{proof}
Fix any $\bm v\in\Hi^{n-1}$. By~\eqref{eq:T-graph} and~\eqref{eq:T-U}, to obtain~\eqref{eq:T-affine-splitting} it suffices to prove
\begin{equation}\label{eq:XC-XU-0}
\X_{\bm C}(\bm v) = \X_{\bm U}(\bm v)+\X_{\bm C}(0).
\end{equation}
We proceed coordinatewise. Using \eqref{eq:XC}, \eqref{eq:traslationPC} and the linearity of $P_{U_i}$, for all $i=1,\ldots,n$, we obtain 
	\begin{align}\label{eq:affine-transf-ind}
		\left(\X_{\bm C}(\bm v)\right)_i
		-\left(\X_{\bm C}(0)\right)_i 
		=~&
		P_{C_i}\left(
		\frac{2}{d_i}\sum_{(h,i)\in\mathcal E}
		\left(\X_{\bm C}(\bm v)\right)_h
		+\frac{1}{d_i}(\ot{Z}\bm v)_i
		\right) \notag \\
		&-P_{C_i}\left(
		\frac{2}{d_i}\sum_{(h,i)\in\mathcal E}
		\left(\X_{\bm C}(0)\right)_h
		\right) \notag \\
		=~&
		P_{U_i}\left(
		\frac{2}{d_i}\sum_{(h,i)\in\mathcal E}
		\left(\X_{\bm C}(\bm v)\right)_h
		+\frac{1}{d_i}(\ot{Z}\bm v)_i
		-\xi_i
		\right)
		+\xi_i \notag \\
		&-P_{U_i}\left(
		\frac{2}{d_i}\sum_{(h,i)\in\mathcal E}
		\left(\X_{\bm C}(0)\right)_h
		-\xi_i
		\right)
		-\xi_i \notag \\
		=~&
		P_{U_i}\left(
		\frac{2}{d_i}\sum_{(h,i)\in\mathcal E}
		\bigl(
		\left(\X_{\bm C}(\bm v)\right)_h
		-\left(\X_{\bm C}(0)\right)_h
		\bigr)
		+\frac{1}{d_i}(\ot{Z}\bm v)_i
		\right).
	\end{align}
Proceeding inductively on the coordinates, one observes that the term~\eqref{eq:affine-transf-ind} coincides with the definition of $(\X_{\bm U}(\bm v))_i$ in~\eqref{eq:XU}.
\end{proof}

The decomposition in \cref{lem:affine-splitting} allows us to study $\T_{\bm C}$ through its linear part $\T_{\bm U}$. We therefore recall the description of the set of fixed points of $\T_{\bm U}$ and its convergence properties established in \cite{graph-linear}. To this end, we first introduce the subspace
\begin{equation}\label{eq: E}
    E:=
    \ot{Z^\dagger}
    \bigl(\mathcal U^\perp\cap\Delta_n^\perp\bigr)
    \subseteq\Hi^{n-1},
\end{equation}
which is closed by \cite[Lemma~4.1]{graph-linear}. Before recalling these results, we establish two equivalent representations of $E$ that will be useful in our analysis.

\begin{lemma}[Alternative representations of $E$]
\label{lem:E-characterization}
The subspace $E$ defined in \cref{eq: E} satisfies
\begin{equation}\label{eq:E-characterization-1}
    E
    =
    \bigl(\ot{Z}^*(\mathcal U)\bigr)^\perp
    =
    \ker\mathcal{X}_{\bm {U}}.
\end{equation}
\end{lemma}

\begin{proof}
Let $\bm e\in\Hi^{n-1}$. By \cite[Remark~4.2]{graph-linear}, we know that $\bm e\in E$ is equivalent to $\ot{Z}\bm e\in\mathcal U^\perp$. Therefore,
\begin{equation*}\label{eq:E-characterization-2}
\begin{aligned}
    \bm e\in E
    &\Longleftrightarrow
    \langle\ot{Z}\bm e,\bm u\rangle=0,
    \quad \forall\bm u\in\mathcal U\\
    &\Longleftrightarrow
    \langle \bm e,\ot{Z}^*\bm u\rangle=0,
    \quad \forall\bm u\in\mathcal U\\
    &\Longleftrightarrow
    \bm e\in\bigl(\ot{Z}^*(\mathcal U)\bigr)^\perp,
\end{aligned}
\end{equation*}
which proves the first equality in~\eqref{eq:E-characterization-1}. In order to prove the second equality, one checks that $\ot{Z}\bm e\in\mathcal U^\perp \Longleftrightarrow P_{\mathcal{U}}(\ot{Z}\bm e)=0 \Longleftrightarrow \mathcal{X}_{\bm{U}}(\bm e)=0$, where the last equivalence follows from~\eqref{eq:XU}. Altogether, this yields $E=\ker\mathcal{X}_{\bm{U}}$.
\end{proof}

\begin{fact}[{\cite[Proposition~4.3 and Theorem~4.4]{graph-linear}}]
\label{th:projectionFixTU}
Let $\T_{\bm U}$ be the homogeneous graph operator given in~\cref{eq:T-U}. Then
\begin{equation*}\label{eq:FixTU}
    \Fix\T_{\bm U}
    =
    \ot{\bm\alpha}(U)\oplus E,
\end{equation*}
where $\bm\alpha\in\mathbb R^{n-1}$ and $E\subseteq\Hi^{n-1}$ are given by \cref{eq:alpha,eq: E}, respectively, and $\oplus$ denotes the Minkowski sum with orthogonal summands. Furthermore,
\begin{equation*}\label{eq:projectionFixTU}
    P_{\Fix\T_{\bm U}}
    =
    \frac{1}{\|\bm\alpha\|^2}
    \ot{\bm\alpha}P_U\ot{\bm\alpha}^*
    +P_E.
\end{equation*}
\end{fact}

\begin{fact}[{\cite[Theorem~4.8]{graph-linear}}]
\label{fact:TU-convergence}
Let $\T_{\bm U}$ be the homogeneous graph operator defined in \cref{eq:T-U}, let $\theta\in{]0,2[}$, and let
$\bm w^0\in\Hi^{n-1}$. Define
\[
    \bm w^{k+1}
    =
    (1-\theta)\bm w^k
    +\theta\T_{\bm U}(\bm w^k),
    \quad \forall k\geq0.
\]
Then \(\bm w^k\to P_{\Fix\T_{\bm U}}(\bm w^0)\).
\end{fact}

\begin{theorem}[Minimal displacement vector in the affine case]
\label{th:affine-minimal-displacement}
Let $C_1,\ldots,C_n$ be closed affine subspaces as in \cref{eq:Ci-affine}, and set $\bm\xi:=(\xi_1,\ldots,\xi_n)$. Then the minimal displacement vector of $\T_{\bm C}$ is
\begin{equation}\label{eq:v-affine}
    v_{\T_{\bm C}}
    =
    P_E\bigl(\ot{Z}^*\bm\xi\bigr),
\end{equation}
where $E\subseteq\Hi^{n-1}$ is given by \cref{eq: E}.
\end{theorem}

\begin{proof}
    Since $\mathcal C=\bm\xi+\mathcal U$, then $\ot{Z}^*(\mathcal C)=\ot{Z}^*\bm\xi+\ot{Z}^*(\mathcal U)$. By \cref{lem:E-characterization} and \cref{thm:graph-ls}, this gives
\[
    v_{\T_{\bm C}}
    =
    P_{\overline{\ot{Z}^*(\mathcal C)}}(0)
    =
    P_{\ot{Z}^*\bm\xi+E^\perp}(0)
    =
    P_E(\ot{Z}^*\bm\xi),
\]
proving the desired equality.
\end{proof}

\begin{remark}
Note that \cref{eq:v-affine} does not depend on the choice of the representatives $\xi_i$ in \cref{eq:Ci-affine}. Indeed, if $\bm\xi'\in\mathcal C$ is another choice (i.e., $\bm\xi'+\mathcal U = \mathcal C = \bm \xi +\mathcal U$), then $\bm\xi'-\bm\xi\in\mathcal U$. Since $E=(\ot{Z}^*(\mathcal U))^\perp$ by \cref{lem:E-characterization}, we get that
\[
    P_E\bigl(\ot{Z}^*\bm\xi'\bigr)
    =
    P_E\bigl(\ot{Z}^*\bm\xi\bigr).
\]
Thus any point of $\mathcal C$ may be used in \cref{eq:v-affine}.
\end{remark}

\begin{remark}[The affine Douglas--Rachford gap vector]\label{rem:DR-affinegap}
In the setting of \cref{rem:DR-range}, suppose that $C_1=\xi_1+U_1$ and $C_2=\xi_2+U_2$. Then $E=U_1^\perp\cap U_2^\perp$ (see~\cite[\S~4.6.1]{graph-linear}), and \cref{th:affine-minimal-displacement} gives
\[
    v_{\T_{\bm C}}
    =
    P_{U_1^\perp\cap U_2^\perp}(\xi_1-\xi_2),
\]
recovering the known representation of the gap vector for two closed affine subspaces; see, e.g.,~\cite[Example~2.2]{BB93}.
\end{remark}

As a direct consequence of \cref{th:affine-minimal-displacement} and \cref{cor:graph-residual}, we obtain the following result.

\begin{corollary}
Let $C_1,\ldots,C_n$ be closed affine subspaces as in \cref{eq:Ci-affine}, and set $\bm\xi:=(\xi_1,\ldots,\xi_n)$. Let $(\bm x^k)_{k\geq1}$ be the shadow sequence generated by \cref{alg:T}. Then
\begin{equation*}\label{eq:affine-increments}
    \ot{Z}^*\bm x^k  \to  P_E\bigl(\ot{Z}^*\bm\xi\bigr),
\end{equation*}
where $E\subseteq\Hi^{n-1}$ is given by \cref{eq: E}.
\end{corollary}

In line with \cref{cor:graph-ls}, the following adds a new characterization of the existence of a solution to the least-squares problem \eqref{eq:graph-ls} that relies on the affine structure of the sets $C_1,\ldots,C_n$.

\begin{corollary}[Existence of least-squares solutions in the affine setting]
\label{prop:affine-ls-existence}
Let $C_1,\ldots,C_n$ be closed affine subspaces as in \cref{eq:Ci-affine}, and set $\bm\xi:=(\xi_1,\ldots,\xi_n)$. Then \eqref{eq:graph-ls} admits a solution if and only if
\[P_{E^\perp}(\ot{Z}^*\bm\xi)\in\ot{Z}^*(\mathcal U).\]
\end{corollary}

\begin{proof}
    By \cref{cor:graph-ls}, it amounts to show that $v_{\T_{\bm C}}\in\ot{Z}^*(\mathcal C)$. Using \cref{th:affine-minimal-displacement} and $\ot{Z}^*(\mathcal C)=\ot{Z}^*\bm \xi + \ot{Z}^*(\mathcal U)$, this is equivalent to
\[
    - P_{E^\perp}(\ot{Z}^*\bm\xi)
    = P_E(\ot{Z}^*\bm\xi)-\ot{Z}^*\bm\xi
    = v_{\T_{\bm C}}-\ot{Z}^*\bm\xi 
    \in\ot{Z}^*(\mathcal U).
\]
Since $\ot{Z}^*(\mathcal U)$ is a linear subspace, the proof is complete.
\end{proof}

\begin{remark}\label{rem:aff-closed}
If $\ot{Z}^*(\mathcal U)$ is closed, then
$\ot{Z}^*(\mathcal U)=E^\perp$ by
\cref{lem:E-characterization}. Hence,
\cref{prop:affine-ls-existence} ensures that
\eqref{eq:graph-ls} admits a solution.
In particular, this holds whenever $\Hi$ is finite-dimensional.
\end{remark}

We now show that, whenever the least-squares problem admits a solution, the shadow sequence converges strongly to a minimizer.

\begin{theorem}[Shadow convergence in the affine case]
\label{th:affine-convergence}
Let $C_1,\ldots,C_n$ be closed affine subspaces as in \cref{eq:Ci-affine}, and suppose that any condition in \cref{cor:graph-ls} holds. Let $(\bm x^k)_{k\geq1}$ be the shadow sequence generated by
\cref{alg:T}. Then
\begin{equation*}\label{eq:shadow-affine-limit2}
    \bm x^k\to\overline{\bm x}
    \in
    \operatorname*{argmin}_{\bm x\in\mathcal C}
    \frac12\|\ot{Z}^*\bm x\|^2.
\end{equation*}
\end{theorem}

\begin{proof}
By \crefpart{cor:graph-ls}{thm:vT-projectionii}, there exists $\bm p\in\Hi^{n-1}$ such that 
\begin{equation}\label{eq:p-normal}
v_{\T_{\bm C}}=\bm p-\T_{\bm C}(\bm p).
\end{equation}
Let $(\bm v^k)_{k\geq0}$ be the governing sequence generated by \cref{alg:T}, and define
\[
    \bm w^k
    :=
    \bm v^k+k\theta v_{\T_{\bm C}}-\bm p,
    \quad \forall k\geq0.
\]
We claim that $(\bm w^k)_{k\geq0}$ is the relaxed fixed-point iteration generated by $\T_{\bm U}$. Indeed, combining \cref{th:affine-minimal-displacement,th:projectionFixTU}, we obtain that $v_{\T_{\bm C}}\in E\subseteq \Fix\T_{\bm U}$. Using this, together with \eqref{eq:T-iteration}, \eqref{eq:p-normal}, \cref{lem:affine-splitting} and the linearity of $\T_{\bm U}$, we obtain
\begin{align*}
    \bm w^{k+1}
    &=
    \bm v^{k+1}+(k+1)\theta v_{\T_{\bm C}}-\bm p\\
    &=
    (1-\theta)\bm v^k
    +\theta\T_{\bm C}(\bm v^k)
    +(k+1)\theta v_{\T_{\bm C}}-\bm p\\
    &=
    (1-\theta)\bm v^k
    +\theta\T_{\bm U}(\bm v^k)
    -\theta\ot{Z}^*\X_{\bm C}(0)
    +(k+1)\theta v_{\T_{\bm C}}-\bm p\\
    &=
    (1-\theta)
    \bigl(\bm w^k+\bm p-k\theta v_{\T_{\bm C}}\bigr)
    +
    \theta\bigl(
        \T_{\bm U}(\bm w^k)+\T_{\bm U}(\bm p)
        -k\theta v_{\T_{\bm C}}
    \bigr)\\
    &\quad
    -\theta\ot{Z}^*\X_{\bm C}(0)
    +(k+1)\theta v_{\T_{\bm C}}-\bm p\\
    &=
    (1-\theta)\bm w^k
    +\theta\T_{\bm U}(\bm w^k)
    +\theta\bigl(
        \T_{\bm U}(\bm p)-\bm p
        -\ot{Z}^*\X_{\bm C}(0)
        +v_{\T_{\bm C}}
    \bigr)\\
    &=
    (1-\theta)\bm w^k
    +\theta\T_{\bm U}(\bm w^k)
    +\theta\bigl(
        \T_{\bm C}(\bm p)-\bm p+v_{\T_{\bm C}}
    \bigr)\\
    &=
    (1-\theta)\bm w^k+\theta\T_{\bm U}(\bm w^k).
\end{align*}
Hence, by \cref{fact:TU-convergence},
\(
    \bm w^k\to\overline{\bm w}
    :=
    P_{\Fix\T_{\bm U}}(\bm v^0-\bm p).
\)

We are now ready to establish the convergence of the shadow sequence.
First, notice that, since $v_{\T_{\bm C}}\in E$, then $\mathcal{X}_{\mathcal{U}}(v_{\T_{\bm C}})=0$ by \cref{lem:E-characterization}.
Using \cref{eq:XC-XU-0} and the linearity of $\X_{\bm U}$, we get
\begin{align*}
    \bm x^{k+1}
    =
    \X_{\bm C}(\bm v^k)
    =
    \X_{\bm U}\bigl(
        \bm w^k+\bm p-k\theta v_{\T_{\bm C}}
    \bigr)
    +\X_{\bm C}(0)
    =
    \X_{\bm U}(\bm w^k+\bm p)+\X_{\bm C}(0).
\end{align*}
Since $\X_{\bm U}$ is continuous, we conclude that
\[
    \bm x^{k+1}
    \to
    \X_{\bm U}(\overline{\bm w}+\bm p)+\X_{\bm C}(0)
    =
    \X_{\bm C}(\overline{\bm w}+\bm p)
    =:\overline{\bm x}\in\mathcal C.
\]
Finally, by \cref{cor:graph-residual} and the continuity of $\ot{Z}^*$, we have $\ot{Z}^*\overline{\bm x}=v_{\T_{\bm C}}$. Thus \cref{thm:graph-ls} shows that $\overline{\bm x}$ solves \cref{eq:graph-ls} and concludes the proof.
\end{proof}

\begin{remark}[Douglas--Rachford shadow convergence]
In the setting of \cref{rem:DR-affinegap}, suppose that a best approximation pair relative to $C_1$ and $C_2$ exists. Then \cref{th:affine-convergence} yields
\[
    \bigl(P_{C_1}(v^k),P_{C_2}\bigl(R_{C_1}(v^k)\bigr)\bigr)
    \to(\overline x_1,\overline x_2)
    \in
    \operatorname*{argmin}_{(x_1,x_2)\in C_1\times C_2}
    \frac12\|x_1-x_2\|^2.
\]
This recovers some of the convergence assertions of \cite[Theorem~4.4]{BM16affine}.
\end{remark}

\section{Illustrative examples}\label{sec:illustrative}

In this section, we  numerically illustrate the results of \cref{sec:affine} on inconsistent feasibility problems formed by three lines in the plane. Thus, we set $\Hi=\mathbb R^2$ and $n=3$, and consider throughout
\begin{equation}\label{eq:Ci-example}
    C_1:= (0,0) + \spa\{(1,0)\}, \;
    C_2:= (0,0) + \spa\{(2,-3)\}, \;
    C_3:= (0,1) + \spa\{(1,-1)\},
\end{equation}
so that $C_1\cap C_2\cap C_3=\varnothing$. We are interested in depicting the minimal displacement vector $v_{\T_{\bm C}}$ for six particular instances of \cref{alg:T}, whose graph settings are gathered in \cref{tab:graph-settings}.  The onto decompositions of $G'$ considered here are those described in~\cite[Fact~2.17]{graph-linear}.

\begin{table}[H]
    \centering
    \renewcommand{\arraystretch}{1.9}\small
    \begin{tabular}{lccc}
        \hline
        \noalign{\vskip -.8ex}Method & $G$ & $G'$ & $Z^*$ \\[-.7ex]
        \hline
        Ryu \cite{ryu20}
        & \agraph{12,13,23}
        & \multirow{2}{*}{\agraph{13,23}}
        & \multirow{2}{*}{\footnotesize $
            \begin{pmatrix}
                1 & 0 & -1 \\[-1ex]
                0 & 1 & -1
            \end{pmatrix}
        $}
        \\[1mm]
        Parallel down \cite{graph-drs}
        & \agraph{13,23}
        &
        &
        \\
        \hline
        Malitsky--Tam \cite{malitsky2023resolvent}
        & \agraph{12,13,23}
        & \multirow{2}{*}{\agraph{12,23}}
        & \multirow{2}{*}{\footnotesize $
            \begin{pmatrix}
                1 & -1 & 0 \\[-1ex]
                0 & 1 & -1
            \end{pmatrix}
        $}
        \\[1mm]
        Sequential \cite{graph-drs}
        & \agraph{12,23}
        &
        &
        \\
        \hline
        \noalign{\vskip .9ex} Parallel up \cite{campoy}
        & \agraph{12,13}
        & \agraph{12,13}
        &\footnotesize   $
            \begin{pmatrix}
                1 & -1 & 0 \\[-1ex]
                1 & 0 & -1
            \end{pmatrix}
        $
        \\[3.5ex]
        \hline
        \noalign{\vskip .9ex} Complete \cite{graph-drs}
        & \agraph{12,13,23}
        & \agraph{12,13,23}
        & \footnotesize $
            \begin{pmatrix}
                \sqrt{2} & -\tfrac{1}{\sqrt{2}} & -\tfrac{1}{\sqrt{2}} \\[-.2ex]
                0 & \sqrt{\tfrac{3}{2}} & -\sqrt{\tfrac{3}{2}}
            \end{pmatrix}
        $
        \\[3.5ex]
        \hline
    \end{tabular}
    \caption{Graph setting defining each splitting algorithm}
    \label{tab:graph-settings}
\end{table}

For each graph setting, \cref{th:affine-minimal-displacement} allows us to compute the minimal displacement vector through the subspace
\[
E=\bigl(\ot Z^*(U_1\times U_2\times U_3)\bigr)^\perp\subseteq (\R^2)^2,
\]
with $U_1=\spa\{(1,0)\}$, $U_2=\spa\{(2,-3)\}$, and $U_3=\spa\{(1,-1)\}$. 
Explicit descriptions of the subspace $E$ for the graph settings considered here were obtained in~\cite[\S~4.5]{graph-linear}. We specialize those formulas to $n=3$ and to the subspaces in~\eqref{eq:Ci-example}. For instance, for Ryu's algorithm, \cite[Eq.~(30)]{graph-linear} gives
\[E = \left\{(e_1,e_2)\in\Hi^2: e_1\in U_1^\perp,\; e_2\in U_2^\perp,\; e_1+e_2\in U_3^\perp\right\}.\]
Since
$U_1^\perp=\spa\{(0,1)\}$,
$U_2^\perp=\spa\{(3,2)\}$, and
$U_3^\perp=\spa\{(1,1)\}$,
we obtain
\[
    E
    =\left\{
      {\bigl(}(0,t),(3s,2s){\bigr)}
      : t,s\in\mathbb R,\ t=s
      \right\}
    =\spa\left\{
      {\bigl(}(0,1),(3,2){\bigr)}
      \right\}.
\]
Moreover, in the present example
\(
\bm\xi={\bigl(}(0,0),(0,0),(0,1){\bigr)},
\)
so that
\(
\ot{Z}^*\bm\xi = {\bigl(}(0,-1),(0,-1){\bigr)}.
\)
Hence \cite[Example~29.18]{bauschke} gives
\[
v_{\T_{\bm C}}
=
P_E(\ot{Z}^*\bm\xi)
=
\frac{
\left\langle
\bigl((0,-1),(0,-1)\bigr),
\bigl((0,1),(3,2)\bigr)
\right\rangle
}{
\left\|
\bigl((0,1),(3,2)\bigr)
\right\|^2
}
\bigl((0,1),(3,2)\bigr)
=
-\frac{3}{14}
\bigl((0,1),(3,2)\bigr).
\]
The corresponding subspaces $E$ and minimal displacement vectors of the remaining methods can be obtained analogously, so we directly collect them in \cref{tab:E-displacement-example}. Note that algorithms sharing the subgraph $G'$ (and using the same $Z$) yield the same subspace $E$. This is the case for Ryu and Parallel down, as well as Malitsky--Tam and Sequential.

\begin{table}[H]
    \centering
    \renewcommand{\arraystretch}{1.1}\small
    \begin{tabular}{c|c|c}
        Method
        & $E$
        & $v_{\T_{\bm C}}=P_E(\ot Z^*\bm\xi)$
        \\
        \hline

        \begin{tabular}{c}
            Ryu\\
            Parallel down
        \end{tabular}
        &
        $\operatorname{span}
        \left\{
        {\bigl(}(0,1),(3,2){\bigr)}
        \right\}$
        &
        $-\frac{3}{14}
        {\bigl(}(0,1),(3,2){\bigr)}$
        \\
        \hline

        \begin{tabular}{c}
            Malitsky--Tam\\
            Sequential
        \end{tabular}
        &
        $\operatorname{span}
        \left\{
        {\bigl(}(0,1),(3,3){\bigr)}
        \right\}$
        &
        $-\frac{3}{19}
        {\bigl(}(0,1),(3,3){\bigr)}$
        \\
        \hline

        \rule[-2.3ex]{0pt}{5.8ex}Parallel up
        &
        $\operatorname{span}
        \left\{
        {\bigl(}(-3,-2),(3,3){\bigr)}
        \right\}$
        &
        $-\frac{3}{31}
        {\bigl(}(-3,-2),(3,3){\bigr)}$
        \\
        \hline

        \rule[-2.3ex]{0pt}{5.8ex}Complete
        &
        $\operatorname{span}
        \left\{
        {\bigl(}(0,\sqrt3),(6,5){\bigr)}
        \right\}$
        &
        $-\frac{3}{32}\sqrt{\frac32}
        {\bigl(}(0,\sqrt3),(6,5){\bigr)}$
    \end{tabular}
    \caption{Subspace $E$ and minimal displacement vector associated with each method for the affine subspaces in~\eqref{eq:Ci-example}}
    \label{tab:E-displacement-example}
\end{table}

We finally turn to the least-squares problems associated with these graph settings. Recall that \eqref{eq:graph-ls} depends only on the subgraph $G'$, and not on the particular onto decomposition of its Laplacian. The resulting least-squares problems for the six algorithms are collected in \cref{tab:graph-methods-ls}. We also include their exact solutions, which can be calculated analytically by parametrizing the affine subspaces and solving the linear system resulting from the first-order optimality condition of the associated unconstrained quadratic problem.

\begin{table}[h]
    \centering
    \renewcommand{\arraystretch}{1.1}
    \resizebox{.98\textwidth}{!}{
    \begin{tabular}{c|c|ccc}
        \multirow{2}{*}{Method}
		&
		\multirow{2}{*}{$
		  \min_{\bm x\in\mathcal C}
		  \frac12\textstyle\sum_{(i,j)\in\mathcal E'}\|x_i-x_j\|^2$}
		&
		\multicolumn{3}{c}{Least-squares solution}
		\\[-1ex]
		&
		&
        $\bar x_1$
        &
        $\bar x_2$
        &
        $\bar x_3$
        \\
        \hline

        \begin{tabular}{c}
            Ryu\\
            Parallel down
        \end{tabular}
        &
        $\displaystyle
        \min_{x_3\in C_3}
        \tfrac12\|x_3-P_{C_1}(x_3)\|^2
        +\tfrac12\|x_3-P_{C_2}(x_3)\|^2$
        &
        { $\left(\tfrac{11}{14},0\right)$}
        &
        { $\left(\tfrac{1}{7},-\tfrac{3}{14}\right)$}
        &
        { $\left(\tfrac{11}{14},\tfrac{3}{14}\right)$}
        \\[1mm]
        \hline

        \begin{tabular}{c}
            Malitsky--Tam\\
            Sequential
        \end{tabular}
        &
        $\displaystyle
        \min_{x_2\in C_2}
        \tfrac12\|x_2-P_{C_1}(x_2)\|^2
        +\tfrac12\|x_2-P_{C_3}(x_2)\|^2$
        &
        { $\left(-\tfrac{2}{19},0\right)$}
        &
        { $\left(-\tfrac{2}{19},\tfrac{3}{19}\right)$}
        &
        { $\left(\tfrac{7}{19},\tfrac{12}{19}\right)$}
        \\[1mm]
        \hline

        \rule[-2.5ex]{0pt}{6.2ex}Parallel up
        &
        $\displaystyle
        \min_{x_1\in C_1}
        \tfrac12\|x_1-P_{C_2}(x_1)\|^2
        +\tfrac12\|x_1-P_{C_3}(x_1)\|^2$
        &
        { $\left(\tfrac{13}{31},0\right)$}
        &
        { $\left(\tfrac{4}{31},-\tfrac{6}{31}\right)$}
        &
        { $\left(\tfrac{22}{31},\tfrac{9}{31}\right)$}
        \\[1mm]
        \hline

        \rule[-2.5ex]{0pt}{6.2ex}Complete
        &
        $\displaystyle
        \min_{\bm x\in\mathcal C}
        \tfrac12 \|x_1-x_2\|^2+\tfrac12\|x_1-x_3\|^2+\tfrac12\|x_2-x_3\|^2$
        &
        { $\left(\tfrac{11}{32},0\right)$}
        &
        { $\left(\tfrac{1}{16},-\tfrac{3}{32}\right)$}
        &
        { $\left(\tfrac{5}{8},\tfrac{3}{8}\right)$}
    \end{tabular}}
    \caption{Least-squares problem~\eqref{eq:graph-ls} for each splitting
    algorithm and its solution for the affine subspaces in~\eqref{eq:Ci-example}}
    \label{tab:graph-methods-ls}
\end{table}

\begin{figure}[htbp]
	\centering

	\begin{subfigure}{\textwidth}
		\centering
		\includegraphics[width=0.75\textwidth]{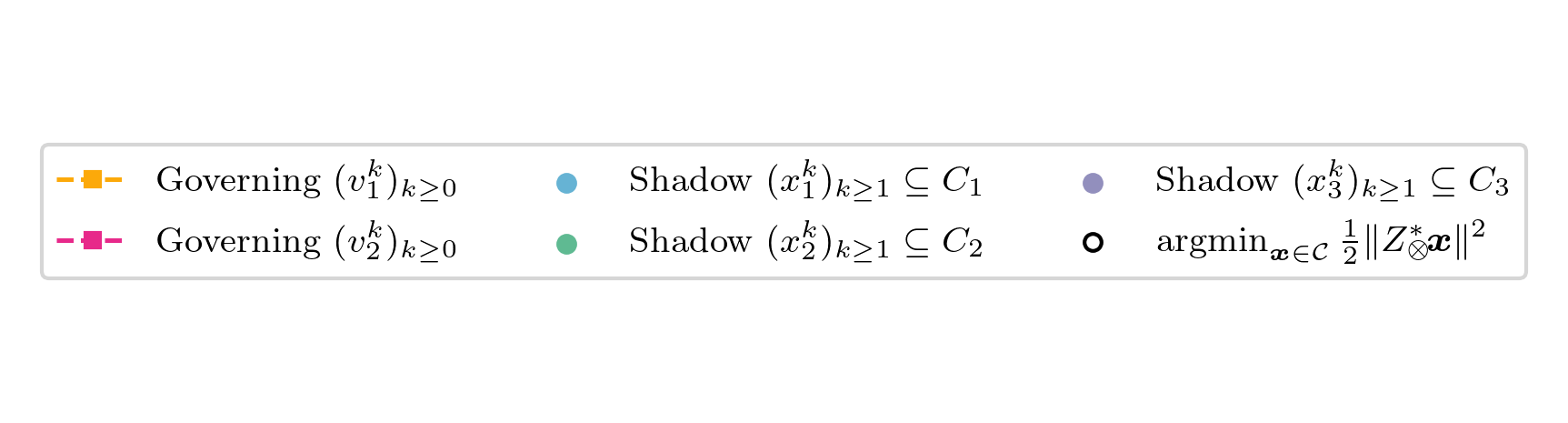} 
	\end{subfigure}
	
	\vspace{-2.5em}
	
	\begin{subfigure}[b]{0.42\textwidth}
		\centering
		\includegraphics[width=\textwidth]{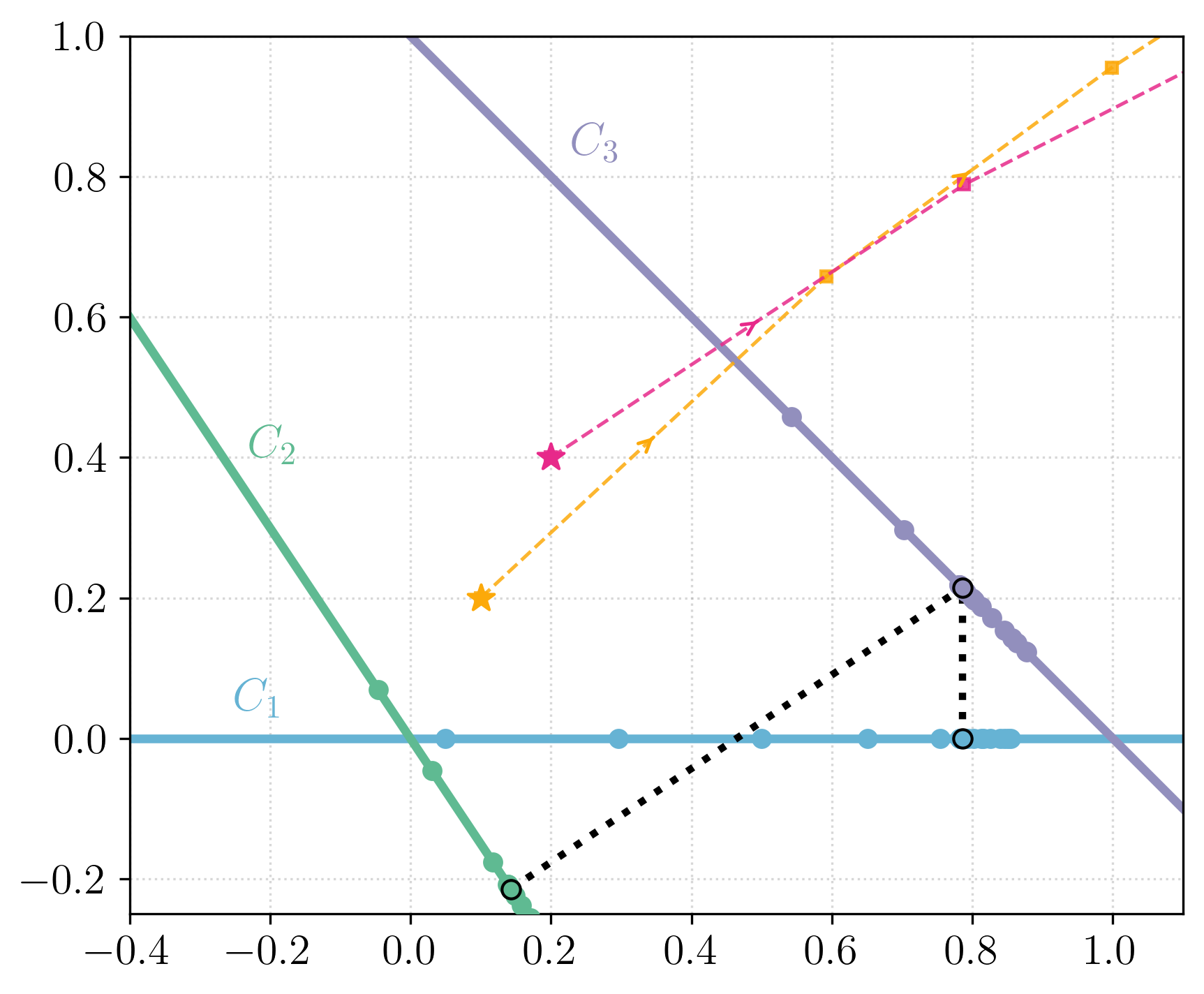}
		\caption{Ryu}
		\label{fig:ryu}
	\end{subfigure}
	\hspace{1mm}
	\begin{subfigure}[b]{0.42\textwidth}
		\centering
		\includegraphics[width=\textwidth]{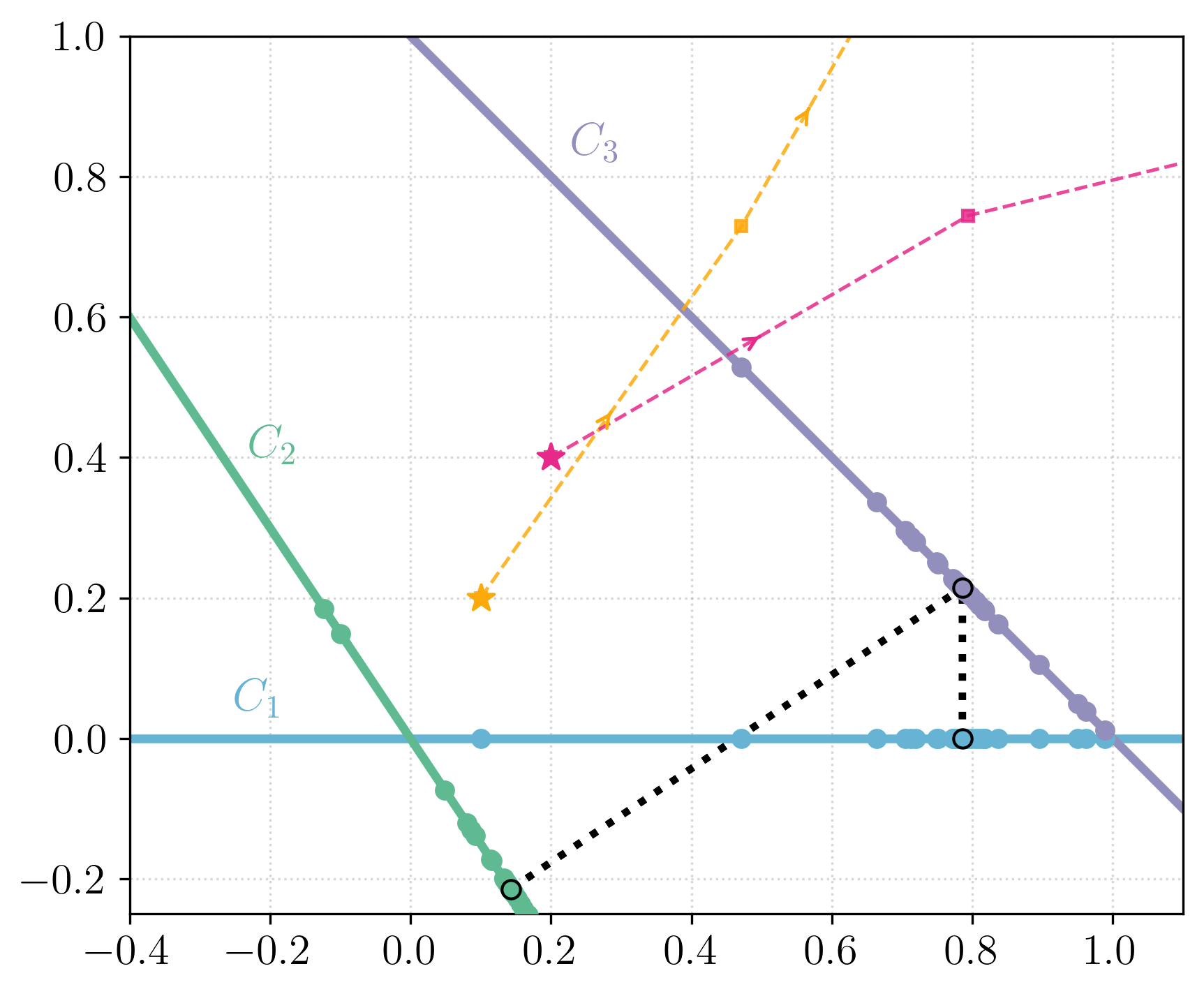}
		\caption{Parallel down}
		\label{fig:paralleldown}
	\end{subfigure}
	
	\vspace{1em}
	\begin{subfigure}[b]{0.42\textwidth}
		\centering
		\includegraphics[width=\textwidth]{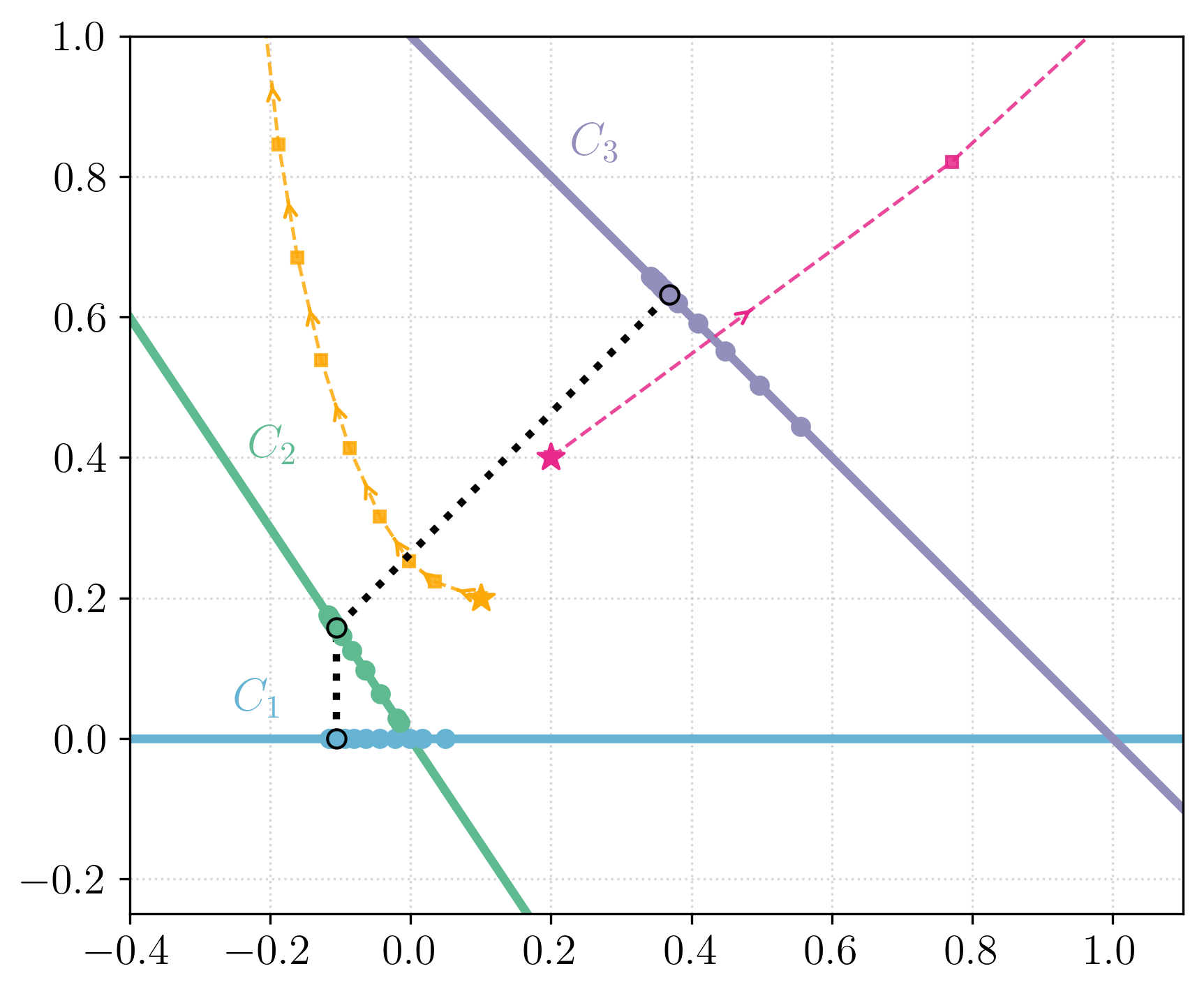}
		\caption{Malitsky--Tam}
		\label{fig:mt}
	\end{subfigure}
	\hspace{1mm}
	\begin{subfigure}[b]{0.42\textwidth}
		\centering
		\includegraphics[width=\textwidth]{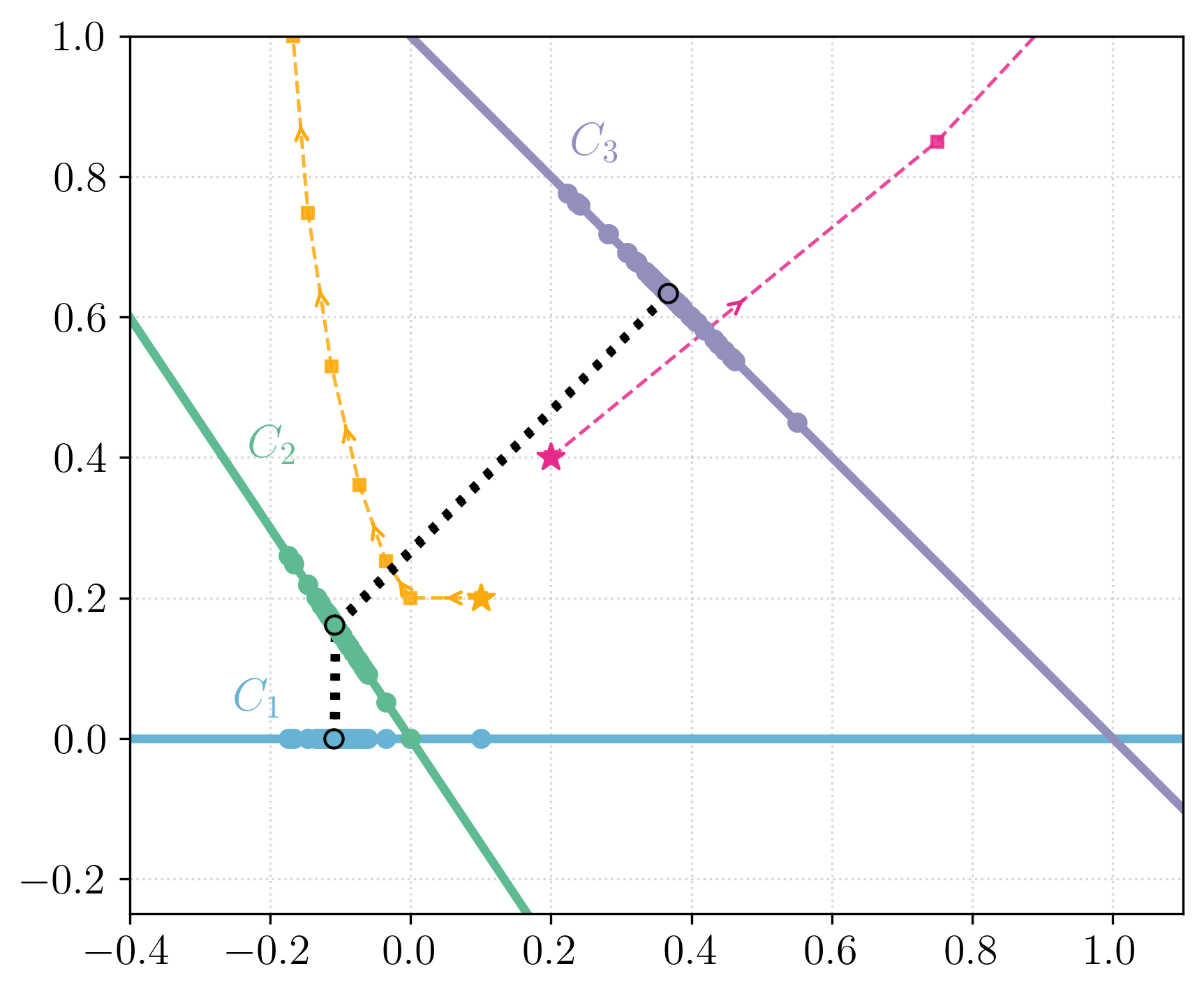}
		\caption{Sequential}
		\label{fig:sequential}
	\end{subfigure}
	
	\vspace{1em}
	\begin{subfigure}[b]{0.42\textwidth}
		\centering
		\includegraphics[width=\textwidth]{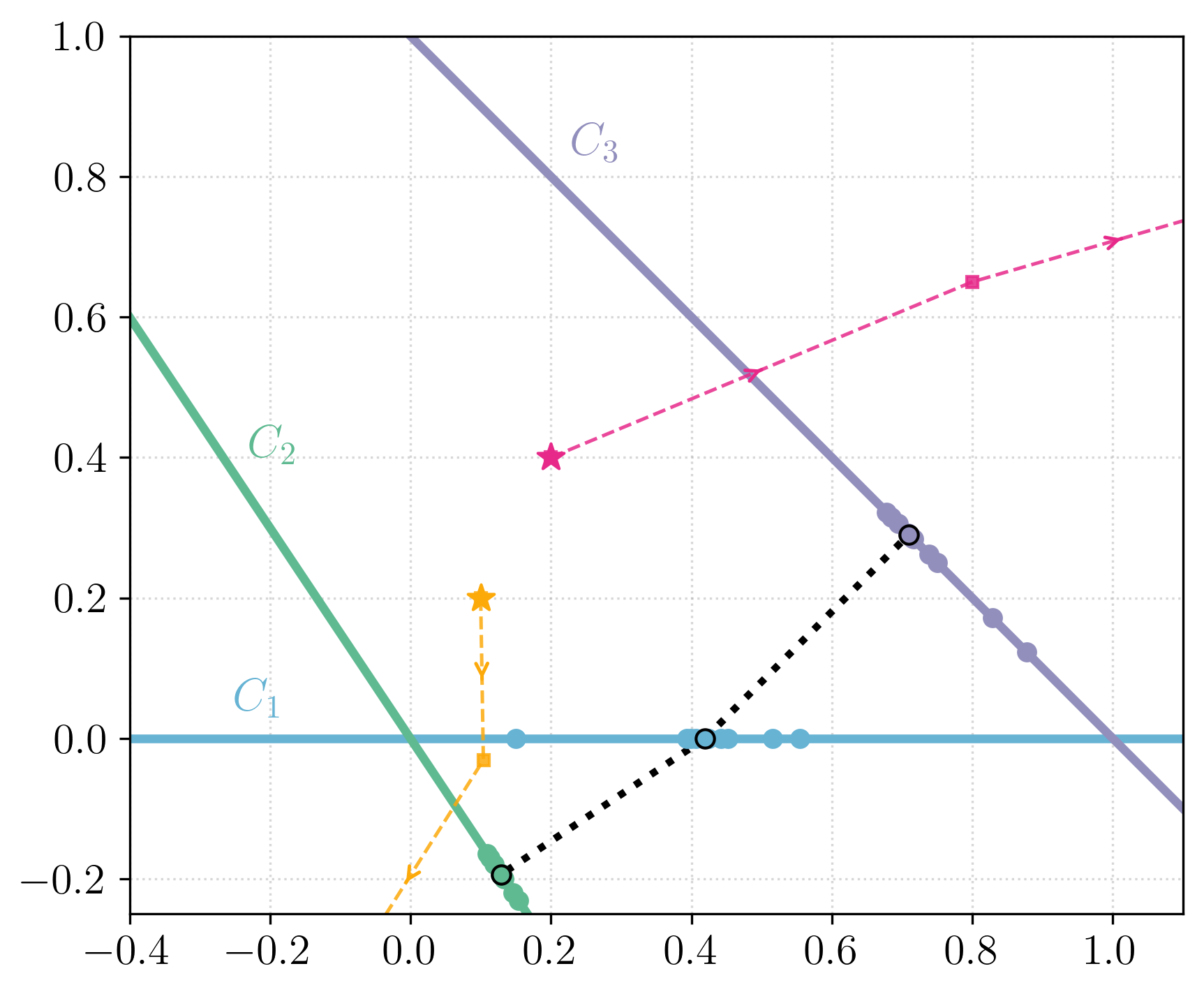}
		\caption{Parallel up}
		\label{fig:parallelup}
	\end{subfigure}
	\hspace{1mm}
	\begin{subfigure}[b]{0.42\textwidth}
		\centering
		\includegraphics[width=\textwidth]{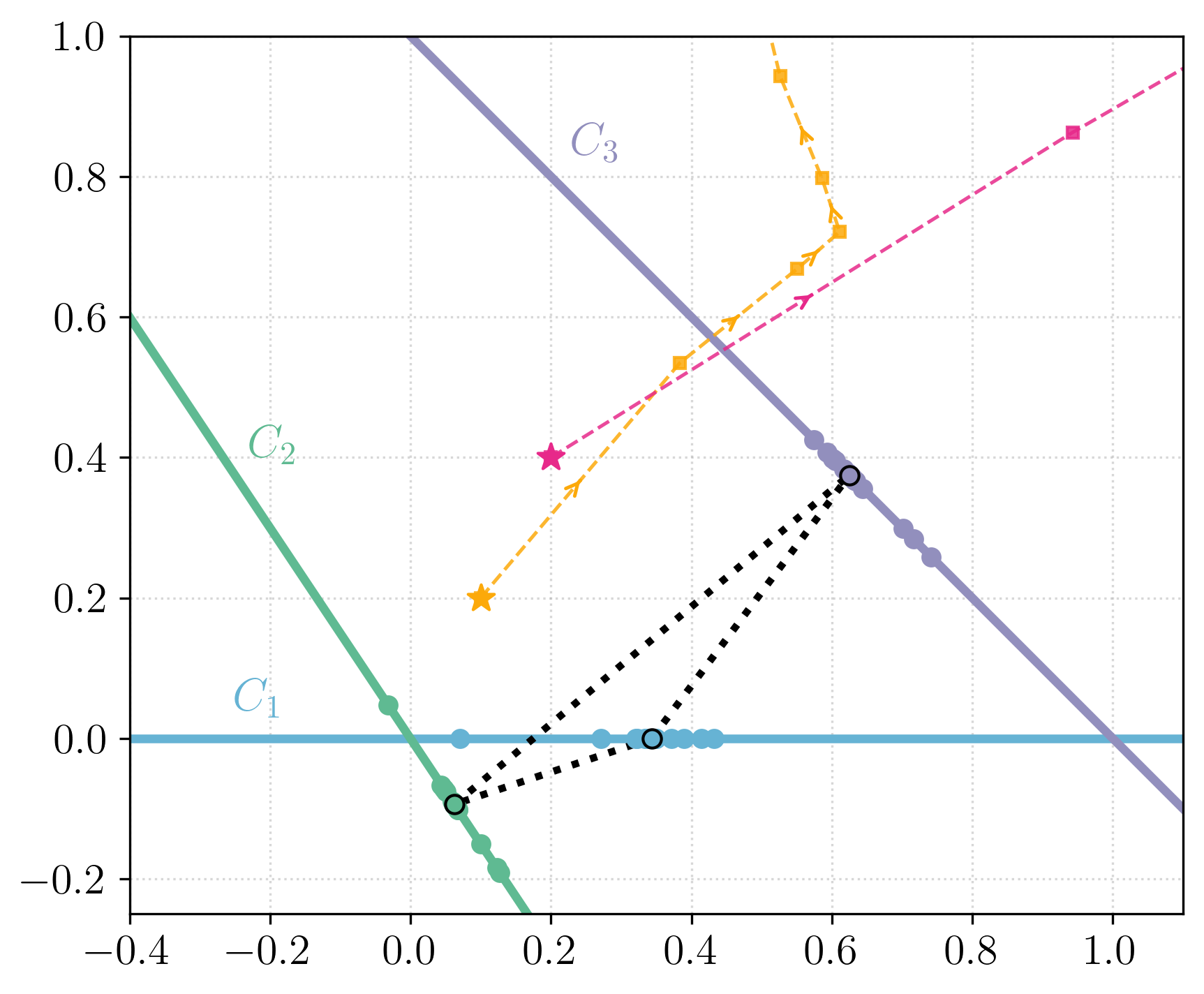}
		\caption{Complete}
		\label{fig:complete}
	\end{subfigure}
	
	\caption{Comparison of three-operator splitting methods for the inconsistent feasibility problem with the affine subspaces in~\eqref{eq:Ci-example}}
	\label{fig:three_operator_comparison}
\end{figure}

By \cref{rem:aff-closed} and \cref{th:affine-convergence}, the shadow sequences of these algorithms converge to a least-squares solution. \cref{fig:three_operator_comparison} illustrates this behavior for each algorithm. In all cases, we use the same relaxation parameter $\theta=1$ and all methods are initialized at the same starting points $v_1^0=(0.1,0.2)$ and $v_2^0=(0.2,0.4)$. We compute up to $50$ iterations of the orbits of the governing sequence $(\bm v^k)_{k\geq 0}$, which are displayed in dashed lines, and of the shadow sequence $(\bm x^k)_{k\geq 1}$, marked as dots. As expected, the shadow sequence converges to the solution of the corresponding least-squares problem, which is represented as an empty-face black dot. Finally, the squared norm of the minimal displacement vector can be obtained as the sum of the squared lengths of the segments drawn in black dotted lines.

\section{Conclusions and future research}\label{sec:conc}
In this work, we studied the behavior of the graph-based splitting framework developed by Bredies, Chenchene and Naldi in~\cite{graph-drs} applied to inconsistent feasibility problems. Thereby, our results extend well-known properties of the Douglas--Rachford splitting algorithm to a wider class of methods, which encompasses popular schemes such as Ryu and Malitsky--Tam, among others. Our main contribution was a characterization of the minimal displacement vector in terms of a least-squares problem measuring disagreement between the sets according to the adjacency of the graph configuration. Furthermore, in the affine setting, we were able to prove strong convergence of the shadow sequence to a least-squares solution, provided it exists. These results were illustrated numerically on several representative graph splitting schemes.
In ongoing work, we aim to investigate whether an analogous weak convergence result can be obtained for general convex sets.

\makeacknowledgements

\bibliographystyle{plain}   
\bibliography{bibliography}

\end{document}